\documentclass[prd,aps,twocolumn,a4paper,floatfix]{revtex4-1}

\usepackage{graphicx,psfrag}
\usepackage{mathrsfs}

\begin{document}

%%%%%%%%%%%%%%%%%%%%%%%%%%%%%%%%%%%%%%%%%%%%%%%%%%%%%%%%%%%%%

\author{Carsten Gundlach}
\affiliation
       {School of Mathematics, 
         University of Southampton \\ 
         Southampton SO17 1BJ, UK \\ 
         Email: cjg@soton.ac.uk}
\author{Jos\'e M.\ Mart\'{\i}n-Garc\'{\i}a}
\affiliation{
        Institut d'Astrophysique de Paris, CNRS,
	 Univ. Pierre et Marie Curie,
	 98bis bd Arago, 75014 Paris, France; \\
	Laboratoire Univers et Th\'eories, CNRS, 
	 Univ. Paris Diderot,
	 5 place Jules Janssen, 92190 Meudon, France \\ 
         Email: jose@xact.es}
\author{David Garfinkle}
\affiliation{Department of Physics, Oakland University, 
         Rochester, MI 48309, USA\\
         and Michigan Center for Theoretical Physics, 
         Randall Laboratory of Physics, University of Michigan, Ann Arbor, MI
         48109-1120, USA \\ 
         Email: garfinkl@oakland.edu}

\title{Summation by parts methods for spherical harmonic
  decompositions of the wave equation in any dimensions}

%%%%%%%%%%%%%%%%%%%%%%%%%%%%%%%%%%%%%%%%%%%%%%%%%%%%%%%%%%%%%

\begin{abstract}

We investigate numerical methods for wave equations in $n+2$ spacetime
dimensions, written in spherical coordinates, decomposed in spherical
harmonics on $S^n$, and finite-differenced in the remaining
coordinates $r$ and $t$. Such an approach is useful when the full
physical problem has spherical symmetry, for perturbation theory about
a spherical background, or in the presence of boundaries with
spherical topology. The key numerical difficulty 
arises from lower-order $1/r$ terms at the origin $r=0$. As a toy
model for this, we consider the flat space linear wave equation in the
form $\dot\pi=\psi'+p\psi/r$, $\dot\psi=\pi'$, where $p=2l+n$, and $l$
is the leading spherical harmonic index. We propose a
class of summation by parts (SBP) finite differencing methods that
conserve a discrete energy up to boundary terms, thus guaranteeing
stability and convergence in the energy norm. We explicitly construct
SBP schemes that are second and fourth-order accurate at interior
points and the symmetry boundary $r=0$, and first and second-order
accurate at the outer boundary $r=R$.

{\it Keywords:} Finite differencing, summation by parts, wave
equation, spherical harmonics.

\end{abstract}

\maketitle

\tableofcontents

%%%%%%%%%%%%%%%%%%%%%%%%%%%%%%%%%%%%%%%%%%%%%%%%%%%%%%%%%%%%%

\section{Introduction}

%%%%%%%%%%%%%%%%%%%%%%%%%%%%%%%%%%%%%%%%%%%%%%%%%%%%%%%%%%%%%

A standard way of proving that the wave equation on flat spacetime
with, for example, Dirichlet boundary conditions is well-posed is to
note that it admits an exactly conserved energy. This energy
functional can then be used to estimate the solution in terms of the
initial and boundary data. The equivalent of well-posedness for the
discretised wave equation is called stability. For suitable
discretisations, stability can be proved in a discrete energy norm
approximating the continuum energy. The Lax equivalence theorem can
then be used to prove convergence in the same norm.

When the background spacetime is curved (as in black hole or stellar
perturbation theory), and/or when the wave equation acquires
lower-order nonlinearities (as in the Einstein equations in
generalised harmonic coordinates), it may still be possible to prove
well-posedness and stability using the existence of a conserved energy
in the constant coefficient approximation to the linearised
equation. See \cite{GKO} for a textbook presentation.

Well-posedness or stability rules out that arbitrarily high frequency
perturbations of the solution grow arbitrarily rapidly. Such
instabilities in finite difference equations appear in practice as
instabilities at the grid frequency that cannot be cured by a small
amount of numerical dissipation. They can, however, be efficiently
eliminated by making sure that the finite
difference scheme conserves a suitable discrete energy when applied to
the linear wave equation in flat spacetime. In the context of
numerical relativity this was shown in a series of papers
\cite{BR1,BR2,BR3}, using finite differencing operators for the wave
equation in Cartesian coordinates proposed by Strand \cite{Strand}.

To show that the time derivative of the energy (integrated over
space) is given only by boundary terms requires integration by
parts. The finite difference operators that preserve a discrete energy
up to boundary terms require an equivalent {\em summation by parts} (from
now, SBP) property. 

In the interior of the numerical domain, Strand's SBP operators are
just the standard symmetric finite-difference operators of minimal
width, for a given order of accuracy. Hence the finite differencing
one would naturally use is already SBP except at the boundaries of the
numerical domain, and in many numerical relativity applications the
outer boundary can be pushed so far out
that problems there can be ignored. This makes it easy to overlook the
importance of the SBP property for stability. By contrast,
\cite{BR1,BR2,BR3,BR4} require full SBP for a clean and stable
treatment of inter-block boundaries in multi-block schemes such as the
``cubed sphere''.

In this paper we develop SBP methods for the wave equation in {\em
  spherical} coordinates. This is natural in three contexts: 1) a
spherically symmetric problem; 2) linear perturbations of a
spherically symmetric background; 3) a physical domain with a
spherical outer boundary. The origin of coordinates then becomes an
unphysical interior boundary $r=0$, which is well-known to cause
numerical instabilities, and which is the major obstacle to using
spherical coordinates. Our methods remove these instabilities at $r=0$
completely and provide a stable treatment of the outer boundary $r=R$.

We do not finite-difference in the angles, but rather start by
decomposing the solution into spherical harmonics. This is natural for
linear equations, where the spherical harmonics decouple, and
can be adapted to the nonlinear case using pseudo-spectral methods.

After the spherical harmonic decomposition and a reduction to first
order (discussed in more detail below) we arrive at the system
\begin{equation}
\label{myeqns}
\dot\psi=\pi', \qquad \dot\pi=\psi'+p{\psi\over r}.
\end{equation}
where the positive integer $p$ is a combination of the dimension of
space and the spherical harmonic index. The finite differencing of
these equations, for $p>0$, is the topic of our paper.

The combination of a spherical harmonic decomposition with finite
differencing in $r$ and $t$ of equations of the type
(\ref{myeqns}) has been used in a number of applications:
spherical gravitational collapse of a scalar field in higher spacetime
dimensions \cite{Birukou,SorkinOren}, gravitational collapse of a
scalar field with angular momentum \cite{Choptuikmwave}, nonspherical
perturbations of spherical relativistic fluid collapse
\cite{fluidpert} and scalar field collapse \cite{critpert}, general
relativistic hydrodynamics \cite{VillainBonazzola}, and Newtonian
magnetohydrodynamics
\cite{TscharnuterWinkler1979,RinconRieutord,IversPhillips}.

The evolution equations (\ref{myeqns}) admit the energy
\begin{equation}
\label{myE}
E\equiv{1\over 2}\int_0^R (\pi^2+\psi^2),\,r^p\,dr
\end{equation}
with time derivative
\begin{equation}
\label{myEdot}
{dE\over dt}=\left(r^p\pi\psi\right)_{r=R} ,
\end{equation}
where Eq.~(\ref{myEdot}) is obtained after using the evolution
equations and the identity
\begin{equation}
\label{integrationbyparts}
\int_a^b\left[\left(\psi'+{p\over
    r}\psi\right)\pi+\pi'\psi\right]r^p\,dr = \left[r^p \pi
  \psi\right]^a_b
\end{equation}
[There is no boundary term at $r=0$ in (\ref{myEdot}) because $\psi$
vanishes there for regular solutions.] The SBP property that our
differential operators need to obey, Eq.~(\ref{SBP}) below, is the
discrete equivalent of (\ref{integrationbyparts}). 

In the linearised Euler equations, for example, (\ref{myeqns}) is
embedded in a larger principal part in such a manner that the identity
(\ref{integrationbyparts}) is still essential for energy
conservation. Hence we believe that SBP operators obeying
(\ref{SBP}) should be used for discretising this piece of the
principal part. However, in the present paper we deal explicitly only
with the wave equation (\ref{myeqns}). 

Underlining our belief that SBP methods are crucial for stability, the
most commonly used second-order accurate discretisation, due to Evans
\cite{Evans}, of the spherical wave equation in 3+1 dimensions (the
case $p=2$), is already SBP in the interior. The SBP approach is {\em
  explicitly} used in \cite{CalabreseNeilsen} to produce a
second-order accurate implementation of the axisymmetric wave equation
(the case $p=1$, see also the work of Sarbach and collaborators
\cite{Sarbach} for a generalisation). Unfortunately, neither of these
methods seems to admit a generalization to higher than second-order
accuracy.  Our contribution is to complete the Evans method to make it
SBP also at the outer boundary $r=R$, to explicitly construct a
fourth-order accurate SBP scheme for any $p$, and to show how schemes
of arbitrary order can be constructed along the same lines.

The plan of the paper is as follows: Sec.~2 presents the
continuum wave equation, the equations that come from its expansion in
spherical harmonics, and our general SBP discretization
framework. Sec.~3 presents our general approach to finding SBP finite
difference operators of arbitrary accuracy, with explicit examples
given of second-order accurate and fourth-order accurate
methods. Sec.~4 treats the outer boundary. Sec.~5 presents numerical
tests of our methods and other methods, while conclusions are
presented in Sec.~6.

%%%%%%%%%%%%%%%%%%%%%%%%%%%%%%%%%%%%%%%%%%%%%%%%%%%%%%%%%%%%%

\section{Continuum equations and their discretization}

%%%%%%%%%%%%%%%%%%%%%%%%%%%%%%%%%%%%%%%%%%%%%%%%%%%%%%%%%%%%%

\subsection{Continuum equations}
\label{subsection:continuum}

In three spatial dimensions, the general solution of the wave equation can
be written in a spherical harmonic series as
\begin{equation}
\Phi(r,t,\theta, \varphi)\equiv \sum_{l=0}^\infty\sum_{m=-l}^l \phi_{lm}(r,t)
Y_{lm}(\theta,\varphi),
\end{equation}
where the partial waves $\phi_{lm}(r,t)$ obey
\begin{equation}
\ddot \phi_{lm} = \phi_{lm}''+{2\over r}\phi_{lm}'-{l(l+1)\over
  r^2}\phi_{lm}.
\end{equation}
A prime denotes $\partial/\partial r$ and a dot $\partial/\partial t$.

This separation of variables ansatz can be generalised to an arbitrary
number of space dimensions.  In polar coordinates, the Laplace
operator $\Delta$ in $n+1$ space dimensions can be split into radial
and angular derivatives as
\begin{equation}
\Delta = \frac{1}{r^n} \frac{\partial}{\partial r}
\left( r^n \frac{\partial}{\partial r}\right) + \frac{1}{r^2} \Delta_{S^n},
\end{equation}
where $\Delta_{S^n}$ is the Laplace operator on the $n$-sphere. (For
what follows we do not need to introduce coordinates on $S^n$.) For
any integer $n\ge 1$, $\Delta_{S^n}$ has eigenfunctions $Y_{l\dots}$ that obey
\begin{equation}
\Delta_{S^n} Y_{l\dots} = -l(l+n-1)Y_{l\dots},
\end{equation}
where $l$ takes integer values $l\ge 0$, and the dots stand for $n-1$
further quantum numbers, for example the index $m$ on $Y_{lm}$ in three
spatial dimensions. We can therefore make the separation of variables ansatz
\begin{equation}
\label{sphericalharmonicexpansion}
\Phi(r,t,{\rm angles})\equiv \sum_{l=0}^\infty\sum_{\dots}
\phi_{l\dots}(r,t) Y_{l\dots}({\rm angles})
\end{equation}
in higher space dimensions, where each partial wave $\phi_{l\dots}$ obeys
\begin{equation}
\label{lwave}
\ddot\phi_{l\dots}=\phi''_{l\dots}+{n\over
  r}\phi_{l\dots}'-{l(l+n-1)\over r^2}\phi_{l\dots}.
\end{equation}
[For $n=1$, (\ref{lwave}) also holds, but $l$ is then the only quantum
  number, is conventionally called $m$, and takes both positive and
  negative integer values.] The restriction to $l=0$, for any $n$,
gives the spherically symmetric wave equation in $n+1$ space
dimensions. From now on, we no longer write the suffix $lm$ or
$l\dots$ that labels the spherical harmonic component $\phi$.

It appears that we have a family of wave equations in $(r,t)$
parameterised by the two integers $n$ (with $n+1$ the dimension of
space) and $l$ (the leading angular quantum number). Considerations of
regularity naturally lead us to an alternative form of this wave
equation in which those two parameters are merged.

We define $\Phi$ to be regular at $r=0$ if and only if it
admits an expansion in positive integer powers of Cartesian
coordinates. When $\Phi$ is expanded in spherical harmonics as in
(\ref{sphericalharmonicexpansion}), this criterion holds if and only
if
\begin{equation}
\label{phibardef}
\phi(r,t)\equiv r^l\bar\phi(r,t), 
\end{equation}
where each $\bar\phi$ admits an expansion in positive {\em even} powers of
$r$. In evolving the wave equation (\ref{lwave}), the condition
$\phi\sim r^l$ is difficult to enforce numerically except for
$l=0,1$. It is easier to evolve $\bar\phi$ itself with the wave
equation
\begin{equation}
\label{lwavebis}
\ddot{\bar\phi}=\bar\phi''+{p\over r}\bar\phi',
\end{equation}
where 
\begin{equation}
p\equiv 2l+n. 
\end{equation}
Recall that $n+1$ is the dimension of space. Hence $p$ is an even
integer in an odd number (in particular, three) of space dimensions,
and an odd integer in an even number of space dimensions. We stress
that in spite of its simple form, this equation represents the wave
equation in any number of spatial dimensions in polar coordinates,
with or without restriction to $SO(n+1)$ symmetry.

The form (\ref{lwavebis}) of our wave equation can further be reduced
to first order in space and time by introducing the auxiliary variables
\begin{equation}
\pi\equiv \dot{\bar\phi}, \qquad \psi\equiv \bar\phi',
\end{equation}
which obey the system (\ref{myeqns}) given in the introduction.
As $\bar\phi$ is an even regular function of $r$, we have
\begin{equation}
\label{evenodd}
\pi(-r,t)=\pi(r,t), \qquad \psi(-r,t)=-\psi(r,t),
\end{equation}
if we formally extend the functions to negative values of $r$.
Generically, $\bar\phi=O(1)$ and hence $\pi=O(1)$ and $\psi=O(r)$ at
the origin.  Eq.~(\ref{myeqns}) is the form of the wave equation that
we will treat for the remainder of the paper, and for which we will
find stable and accurate finite difference numerical approximations.

In order to control the growth of $E$, the boundary term at $r=R$ must
be controlled by a suitable boundary condition. Here we consider outer
boundary conditions of one of three forms. (For simplicity, we
  consider only homogeneous boundary conditions.) The
well-known maximally dissipative boundary conditions are
\begin{equation}
\label{maxdissBC}
\rho\pi+\sigma\psi=0, \quad r=R, \quad \rho\sigma\ge 0.
\end{equation}
From (\ref{myEdot}) it is clear that these give $dE/dt\le 0$.
We also consider the higher-order boundary conditions
\begin{equation}
\label{piBC}
\rho\pi+\mu\pi'=0, \quad r=R, \quad \rho\mu\ge 0,
\end{equation}
or
\begin{equation}
\label{psiBC}
\sigma\psi+\nu\left(\psi'+{p\over r}\psi\right)=0, \quad r=R, \quad
\sigma\nu\ge 0.
\end{equation}
Appendix~\ref{appendix:BCderivslwave} shows that these make a
  modified energy nonincreasing. Hence the wave equation with any of
these boundary conditions is well-posed. A continuum energy exists and
implies well-posedness also for the more general class of boundary
conditions
\begin{equation}
\rho\pi+\sigma\psi+\mu\pi'+\nu(\psi'+p\psi/r)=0, 
\end{equation}
for certain parameter ranges, but we have not been able to find a
discrete counterpart for this case. 

%%%%%%%%%%%%%%%%%%%%%%%%%%%%%%%%%%%%%%%%%%%%%%%%%%%%%%%%%%%%%

\subsection{Discretisation}
\label{section:discretisation}

Throughout this paper we finite-difference in $r$ only, but
assume the continuum limit in time. A fully discrete scheme is
obtained at the end by using a suitable ODE solver in $t$ (the method
of lines).

We use grid functions $\Psi_i(t)$ and $\Pi_i(t)$ on a grid $r_i$ to
represent the continuum functions $\pi(r,t)$ and $\psi(r,t)$, assuming
that $\Pi_i(t)\equiv\pi(r_i,t)$ and $\Psi_i(t)\equiv\psi(r_i,t)$, and
that $\pi(r,t)$ and $\psi(r,t)$ admit Taylor expansions in $r$ to the
required order at any $r$.  From now on, we suppress the
$t$-dependence as it is relevant only later when we add time
discretisation using the method of lines, that is we write $\pi(r)$
and $\Pi_i$, etc. We also use a matrix notation where grid
functions are written as column vectors, e.g. $\Pi$, and finite
differencing operators as matrices acting on these vectors,
e.g. $D\Pi$. 

A $(2K+1)$-point difference operator $\tilde D$ is defined by
\begin{equation}
(\tilde D \Psi)_i=\sum_{j=i+s-K}^{i+s+K} \tilde D_{ij} \Psi_j
\end{equation}
where $-K\le s\le K$ is an offset. The parameters $\tilde D_{ij}$ of
the difference operator are simply the elements of the band-diagonal
matrix $\tilde D$. 

We assume a uniform grid with step size $\Delta r\equiv h$. 
Our methods will require a grid that is either staggered 
or centred about $r=0$. In either case
we find it convenient to introduce the notation
\begin{equation}
r_i \equiv ih, \qquad i={1\over 2},{3\over 2},\dots,M \quad \hbox{or}\quad
 i=0,1,\dots M,
\end{equation}
that is, the grid index $i$ takes half-integer values for the
staggered grid and integer values for the centred grid. In either case
$R\equiv r_M\equiv Mh$.  Whenever needed, we formally extend the grid functions
to any negative value of $i$ with $\Psi_{-i}=-{\Psi_i}$ and
${\Pi_{-i}}={\Pi_i}$.

%%%%%%%%%%%%%%%%%%%%%%%%%%%%%%%%%%%%%%%%%%%%%%%%%%%%%%%%%%%%%

\subsection{Summation by parts}
\label{subsection:SBP}

As is well-known, the continuum equations (\ref{myeqns}) are
well-posed in the norm provided by $E$, given in (\ref{myE}) above,
because $E$ is conserved. A summation by parts (SBP) finite
differencing scheme exactly conserves a discrete equivalent $\hat E$
of the continuum energy $E$. This guarantees that it is stable (the
discrete equivalent of well-posed) in the energy norm.

We consider the discrete energy
\begin{equation}
\label{hatEdef}
{\hat E}\equiv{1\over 2}h^{p+1}\left(\Pi^t W\Pi+\Psi^t \tilde W\Psi\right),
\end{equation}
where ${}^t$ denotes the matrix transpose and where
\begin{equation}
\tilde W^t = \tilde W, \quad W^t=W, \quad \tilde W>0, \quad W>0,
\end{equation}
and we write the finite differencing scheme as
\begin{equation}
\dot \Psi=h^{-1}D\Pi, \quad \dot\Pi=h^{-1}\tilde D\Psi.
\end{equation}
The powers of $h$ have been introduced so that $W$, $\tilde W$, $D$,
$\tilde D$ are all dimensionless and independent of $h$. The quantity
$r_i/h=i$ also has this property. We will derive explicit
expressions later, but both $h^pW$ and $h^p\tilde W$ approximate
$r^p$, while $h^{-1}D$ approximates $d/dr$ and $h^{-1}\tilde D$
approximates $d/dr+p/r$.

The SBP property that guarantees that
$\hat E$ is constant up to boundary terms is
\begin{equation}
\label{SBP}
W\tilde D+(\tilde WD)^t=B,
\end{equation}
where the boundary operator $B$ is defined by
\begin{equation}
\label{Bdef}
\Pi^tB\Psi \equiv \chi M^p \Pi_M\Psi_{M},
\end{equation}
and the constant $\chi$ obeys $\chi\to 1$ in the continuum limit
$M\to\infty$ as $h\to 0$ at fixed $r=R$. [There is no boundary
  contribution at $r=0$, consistent with the fact that we impose
  $\psi(0)=0$.]  Eq. (\ref{SBP}) is the discrete equivalent of
(\ref{integrationbyparts}).

As $W$ is positive definite, it is invertible, and we can consider
$\tilde D$ as determined by a choice of $D$, $W$, $\tilde W$ and $B$:
\begin{equation}
\label{SBPbis}
\tilde D=- W^{-1}D^t \tilde W+W^{-1}B.
\end{equation}

In the case $p=0$ considered by Strand \cite{Strand}, $W$ and $\tilde
W$ represent $1$, and $\tilde D$ and $D$ both represent $d/dr$. It is
then natural to set $W=\tilde W$ and $D=\tilde D$. 

%%%%%%%%%%%%%%%%%%%%%%%%%%%%%%%%%%%%%%%%%%%%%%%%%%%%%%%%%%%%%%%%%%%%%%%%%%

\subsection{The symmetry boundary $r=0$} 

In numerical simulations using polar coordinates one is faced with the
fact that $r=0$ is a boundary of the numerical grid, but is not in
fact a boundary of the physical domain. As a result, there are
(typically) no physical boundary conditions one can or must impose in
the continuum limit, but the numerical simulation does require
boundary conditions. These are derived from the assumption that the
desired solution is not less differentiable at $r=0$ than for
$r>0$. As stated earlier, we assume $\Phi$ to be smooth in Cartesian
spatial coordinates, which is equivalent to $\pi$
being smooth and even and $\psi$ being smooth and odd.  The standard
general approach to imposing such ``symmetry boundary conditions'' or
``regularity conditions'' is to extend the numerical grid into a small
number of ``ghost points'' representing negative $r$ which are
populated by the assumed even or odd parity of the grid
functions. Standard centred finite differencing methods can then be
used at and near the boundary as if it was an interior point.

From a strict SBP point of view, there are no ghost points, and finite
difference operators are necessarily skewed near the boundary. The
fact that $r=0$ is not a physical boundary is represented by the fact
that $B$ is zero at the boundary $r=0$. 

However, we find that the use of ghost points as a notational device
allows a simpler derivation, presentation, and application of our
results, in that we do not need to discuss $r=0$ explicitly as a
boundary. Rather than introduce a few ghostpoints, for our {\em
  derivation} we extend all grid objects from $1/2$ or $0,\dots,M$ to
$-M,\dots M$, corresponding to $-R\le r\le R$. We can then formally
treat $r=0$ as an interior point.

We extend the grid functions to negative $i$ as
\begin{equation}
\label{discreteevenodd}
\Pi_{-i}=\Pi_i, \qquad \Psi_{-i}=-\Psi_i.
\end{equation}
Because $W$ and $\tilde W$ are used only to define $\hat E$, we can assume
without loss of generality that
\begin{equation}
\label{Wsymmetry}
W_{-i,-j}\equiv W_{ij}, \qquad W_{-i,j}\equiv W_{i,-j}=0, \qquad i,j>0,
\end{equation}
and similarly for $\tilde W$. $B$ is extended by $B_{-M,-M}=-B_{MM}$.
In Appendix~\ref{appendix:ghostpoints} we prove from these assumptions
that (\ref{discreteevenodd}) holds at all times if and only if
\begin{equation}
\label{Dsymmetry}
D_{-i,-j}=-D_{ij}.
\end{equation}

When {\em coding} our method, we implement $D$ and $\tilde D$ with a
few ghost points. Equivalently, the ghost points can be explicitly
eliminated. A rigorous discussion of this point is relegated to
Appendix~\ref{appendix:ghostpoints}, as it introduces additional
notation not required for our main argument. Obviously, our time
updates will by construction exactly preserve the evenness of $\pi$ and
oddness of $\psi$.

%%%%%%%%%%%%%%%%%%%%%%%%%%%%%%%%%%%%%%%%%%%%%%%%%%%%%%%%%%%%%

\section{Accuracy}
\label{section:accuracy}

%%%%%%%%%%%%%%%%%%%%%%%%%%%%%%%%%%%%%%%%%%%%%%%%%%%%%%%%%%%%%

\subsection{General considerations}

In this section, we will consider only the behavior of the finite
difference operators at interior points (including $r=0$) of the
numerical grid, postponing to the next section the discussion of how
the operators behave at and near the outer boundary. In what follows,
we will always choose the finite difference operator $D$ to be a
standard centred difference operator of the appropriate order. That is,
for second-order accurate methods, for interior points, we will choose
\begin{equation}
{{(D\Pi})_i} \equiv {\frac {{\Pi _{i+1}} - {\Pi _{i-1}}} 2},
\label{D2}
\end{equation}
while for fourth-order accurate methods, we will choose
\begin{equation}
\label{D4}
(D\Pi)_i \equiv {\frac {8(\Pi_{i+1}-\Pi_{i-1})-(\Pi_{i+2}-\Pi_{i-2})} {12}}.
\end{equation}
Once we choose $W$ and $\tilde W$, the operator $\tilde D$ is given by
Eq.~(\ref{SBPbis}) and the scheme preserves the discrete energy of
Eq.~(\ref{hatEdef}), and thus is stable.  Our task then is to choose
$W$ and $\tilde W$ in such a way that the operator ${h^{-1}} {\tilde
  D}$ so determined is an accurate (to the chosen order) finite
difference representation of the continuum operator $d/dr + p/r$.  

In analyzing the accuracy of $\tilde D$ it is helpful to write the
grid values $\Psi_j=\psi(r_j)$ in terms of the Taylor expansion of
$\psi(r)$ about the fixed grid point $r_i$. We can then write
\begin{eqnarray}
\label{cdef}
h^{-1}(\tilde D\Psi)_i=c_{0i} h^{-1}
\psi(r_i)+c_{1i}\psi'(r_i)+c_{2i} h\psi''(r_i) \nonumber \\ +\dots +
h^{2K-1}c_{2K,i}\psi^{(2K)}(r_i)+O\left(h^{2K}\right), \quad
\end{eqnarray}
where the $c_{\alpha i}$ are a set of numbers linearly related in a
straightforward way to the $\tilde D_{ij}$ at each point $i$.  In the
following we adopt a simplified notation where the $c_{\alpha i}$
(with $\alpha=0,\dots,2K$) are written as $c_\alpha$, $\psi(r_i)$
simply as $\psi$, etc., and $r_i$ simply as $r$. That is, we do not
write the dependence on $i$, and all continuum quantities are
evaluated at $r=r_i$.

The difference operator $\tilde D$ is said to be accurate to order
$2N$ if it obeys (using our abbreviated notation)
\begin{equation}
\label{psiaccuracy}
h^{-1}(\tilde D\Psi)_i = {p\over r}\psi+\psi'+O\left(h^{2N}\right). 
\end{equation}
The point at $r=0$, which arises (only) on a centred grid, must be
treated specially.  Taking the limit as $r \to 0$ of
Eq.~(\ref{psiaccuracy}) at finite $h$ we see that at $r=0$
\begin{equation}
{c_1} = 1 + p
\label{cequalities_origin}  
\end{equation}
while the other odd $c_\alpha$ vanish and the even $c_\alpha$ are undetermined.

A key observation for what follows is that (\ref{psiaccuracy}) 
needs to be obtained formally in the limit $h\to 0$, {\em both} at
(approximately) constant $r$, and at constant $i$. The possible
problem with the latter limit are error terms of the form $h^m/r^n$,
which are $O(h^m)$ at constant $r$, but only $O(h^{m-n})$ at constant
$i$.  

Naively one would expect the accuracy requirement (\ref{psiaccuracy})
for $\tilde D$ at $r\ne 0$ to be equivalent to the following
constraints on the coefficients of the difference operator (as defined
above):
\begin{equation}
\label{cequalities}
c_0={ph\over r}, \quad c_1=1, \quad c_2=\dots=c_{2N}=0.
\end{equation}
Clearly, we would need a stencil of width $2N+1$ or larger to control
all these $c_\alpha$, as, for $p>0$, the even $c_\alpha$ cannot be set
to zero just by using a symmetric stencil. However, we shall now see
that we can violate some of the equalities (\ref{cequalities}) as
$r\to 0$ and in effect replace them with approximate equalities. The
effect is that we will only need an $N+1$ point stencil.

Rather than devising a general notation, we present the cases $N=1$
and $N=2$, after which it should be clear how one can proceed to
arbitrary $N$. 

For $N=1$, we make the following ansatz:
\begin{eqnarray}
\label{N1first}
c_0 &=& {ph\over r}+\delta_0\left({h\over r}\right)^3, \\
c_1 &=& 1-\delta_0\left({h\over r}\right)^2, \\
\label{N1last}
c_2 &=& \delta_1\left({h\over r}\right),
\end{eqnarray}
where the $\delta_\alpha$ may depend on $i$.  The special case
$\delta_0=\delta_1=0$ brings us back to (\ref{cequalities}), but we
shall now see that the parameters $\delta_\alpha$ do not need to
vanish identically but only need to be bounded because of the way
$\psi'$ approximates $\psi/r$ and vice versa for regular odd functions
$\psi(r)$ as $r\to 0$.  Substituting this ansatz into (\ref{cdef})
gives
\begin{eqnarray}
\label{ansatz2}
h^{-1}(\tilde D\Psi)_i&=&{p\over r}\psi+\psi' 
+\delta_0 h^2 \left[r^{-2}\left({\psi\over r}-\psi'\right)\right]\nonumber \\
&& +\delta_1 h^2 \left[r^{-1}\psi''\right]+R_2.
\end{eqnarray}
Here 
\begin{equation}
\label{R2def}
R_2=c_3h^2\psi'''+c_4h^3\psi''''+\dots,
\end{equation}
where for a 3-point stencil $c_3$, $c_4$, \ldots are known linear
functions of $c_0$, $c_1$ and $c_2$.  Now, because $\psi$ can be
expanded in positive odd integer powers of $r$, both square brackets
in (\ref{ansatz2}) are actually $O(1)$ as $r\to 0$.  Therefore, as
long as $\delta_0$ and $\delta_1$ are bounded uniformly in $i$, the
coefficients of $h^2$ in (\ref{ansatz2}) are bounded uniformly in
$i$. Similarly, as $c_3,c_4,\dots$ are regular functions of $\delta_0$
and $\delta_1$, the coefficients of $h^2$ and all higher powers of $h$
in (\ref{R2def}) are also explicitly regular at $r=0$ and so we have
the desired second-order accuracy, uniformly in $i$.

For $N=2$ we make the ansatz
\begin{eqnarray}
\label{N2first}
c_0 &=& {ph\over r}+\delta_0\left({h\over r}\right)^5, \\
c_1 &=& 1-\delta_0\left({h\over r}\right)^4, \\
c_2 &=& \left({\delta_0\over 3}+\delta_1\right)\left({h\over
  r}\right)^3, \\
c_3 &=& -\delta_1\left({h\over r}\right)^2, \\
\label{N2last}
c_4 &=& \delta_2\left({h\over r}\right),
\end{eqnarray}
which gives
\begin{eqnarray}
\label{ansatz4}
h^{-1}(\tilde D\Psi)_i&=&{p\over r}\psi+\psi' \nonumber \\
&& +\delta_0 h^4
\left[r^{-3}\left({\psi\over r^2}-{\psi'\over r}+{\psi''\over
    3}\right)\right]\nonumber \nonumber \\
&&+\delta_1 h^4
\left[r^{-2}\left({\psi''\over r}-\psi'''\right)\right]\nonumber \\ &&
+\delta_2 h^4 \left[r^{-1}\psi ''''\right]+R_4,
\end{eqnarray}
where $R_4=O(h^4)$ in the sense discussed above.
Again, all the square brackets are regular at $r=0$, and so we have
fourth-order accuracy if and only if the $\delta_{\alpha i}$ are bounded
uniformly in $i$.

It should now be clear that this method can be extended to arbitrary
$N$, giving $N$ equations to be solved through a suitable choice of $W$ and
$\tilde W$, and $N+1$ inequalities (uniform in $i$ bounds on the
$\delta_{\alpha i}$) to be then verified for that solution.

Informally, our method can be described as ``trading $r$ for $h$''. It
works because the terms in square brackets above are all $O(1)$ as
$r\to 0$, which in turn requires $\psi(r)$ to be a regular odd
function of $r$.

Our task has thus become to choose $W$ and $\tilde W$ in such a way
that the operator $\tilde D$ given by equation (\ref{SBPbis})
satisfies our ansatz [Eqs.~(\ref{N1first}-\ref{N1last}) for $N=1$ and
Eqs.~(\ref{N2first}-\ref{N2last}) for $N=2$] such that the
quantities $\delta_\alpha$ are uniformly bounded. We now show explicitly
how this task can be accomplished.

%%%%%%%%%%%%%%%%%%%%%%%%%%%%%%%%%%%%%%%%%%%%%%%%%%%%%%%%%%%%%

\subsection{Second-order accuracy (SBP2)}

We begin with the case $N=1$.  For simplicity, we choose $W$ and
$\tilde W$ to be diagonal. That is,
\begin{equation}
\label{vwdef}
W={\rm diag}(w_i), \qquad \tilde W={\rm diag}(v_i).
\end{equation}
The SBP formula (\ref{SBPbis}) then gives
\begin{equation}
\label{Dtilde2}
(\tilde D\Psi)_i={v_{i+1}\Psi_{i+1}-v_{i-1}\Psi_{i-1}\over 2w_i}
\end{equation}
for interior points. 

We have allowed for $v_i\ne w_i$ because this allows us to cover the
Evans and Sarbach methods reviewed in the Appendix, but for
the remainder of this Subsection we further restrict our ansatz to $v_i=w_i$,
using $w_i$ as the parameters.
We can then read off $c_0$, $c_1$ and $c_2$ in terms of $w_i$. 
The one equality contained in (\ref{N1first}-\ref{N1last}), namely
\begin{equation}
\label{1+p}
\left({r\over h}\right)c_0+c_1=1+p,
\end{equation}
keeping in mind that $r/h=i$, gives a linear recurrence relation of
degree 2 for $w_i$,
\begin{equation}
\label{SBP2recurrence}
(i+1)w_{i+1}-(i-1)w_{i-1}=2(p+1)w_i.
\end{equation}
The other two accuracy conditions define $\delta_1$ and $\delta_2$ in
terms of $w_i$. On a {\em staggered} grid, from (\ref{Wsymmetry}) we
have $w_{-1/2}=w_{1/2}$. We initially fix an arbitary value for
$w_{1/2}$, and can then solve the recursion for $w_i$ for all $i\ge
3/2$. (Note that $\tilde D$ is unchanged if $W$ and $\tilde W$ are
multiplied by the same constant factor). On a {\em centred grid},
evaluating Eq.~(\ref{cequalities_origin}) with $w_{-1}=w_1$ gives
$w_1=(1+p)w_0$. We initially fix an arbitrary value of $w_0$ and can
then solve the recursion for $w_i$ for all $i\ge 2$.

The $w_i$ determine the operator $\tilde D$ which in turn determines
the quantities $\delta _0$ and $\delta _1$. These quantities are
plotted in Fig.~\ref{deltas_N=1}. Note that these quantities are
uniformly bounded, which confirms that our method is second-order
accurate.  Appendix~\ref{appendix:recurrence} confirms this
analytically. For comparison, Fig.~\ref{deltas_N=1} also contains the
corresponding quantities for the method of Evans \cite{Evans}, which we
present in our notation in Appendix~\ref{appendix:Evans}.

%%%%%%%%%%%%%%%%%%%%%%%%%%%%%%%%%%%%%%%%%%%%%%%%%%%%%%%%%%%%%

\subsection{Fourth-order accuracy (SBP4)}

We now turn to the case of $N=2$, that is a fourth-order accurate scheme.
We can no longer choose $W$ and $\tilde W$ to be identical and diagonal.
Instead, we choose $W$ to be diagonal and $\tilde W$ to be band-diagonal with three
bands. We parameterize them as 
\begin{eqnarray}
W_{i,i}&=&w_i, \qquad w_{-i}=w_i, \\
\tilde W_{i,i}&=& v_i, \qquad v_{-i}=v_i, \\
\tilde W_{i,i+1}&=& u_{i+1/2}, \qquad u_{-i}=u_i, \\
\tilde W_{i,i-1}&=& u_{i-1/2},  
\end{eqnarray}
and all other components zero, where on the staggered grid the
index on $v$ and $w$ takes half-integer values and the index on $u$
takes integer values, and the other way around on the centred grid.  
In the interest of simplicity, we would like to have as few nonvanishing
$u_i$ as possible.  On the staggered grid it is possible to have only 
$u_1$ nonvanishing, while on the centred grid, it is possible to make only
$u_{3/2}$ and $u_{5/2}$ nonvanishing.  From now on, we make this choice
of $u_i$.  

The ansatz of Eqs.~(\ref{N2first}-\ref{N2last}) imply two equalities,
namely Eq.~(\ref{1+p}) and
\begin{equation}
{c_1} + 3 {\frac h r} {c_2} + 3 {{\left ( {\frac h r} \right ) }^2} {c_3} = 1.
\label{equality2}
\end{equation}
If we temporarily take $u_i$ as given, Eqs.~(\ref{1+p}) and
(\ref{equality2}) determine the $w_i$ plus a linear recurrence
relation of order 4 for the $v_i$. On the {\em staggered} grid, we can fix
$v_{1/2}=v_{-1/2}$ and $v_{3/2}=v_{-3/2}$ arbitrarily, and solve the
recurrence relation for $v_i$ for $i\ge 5/2$ starting from those four
points and our choice of $u_1$.  On the {\em centred} grid, the accuracy
conditions (\ref{cequalities_origin}) at the origin reduce to
$(1+p)w_0=v_1-(1/8)u_{3/2}+(5/8)u_{5/2}$ and
$v_2=v_1+(63/8)u_{3/2}-(27/8)u_{5/2}$.  We can fix $v_1$ and choose
$u_{3/2}$ and $u_{5/2}$ arbitrarily and then compute $v_i$ for $i\ge
3$ from the recurrence relation. (Note that $v_0$ multiplies $\Psi_0$,
which vanishes, and hence does not participate in the recurrence.) It
remains to fix the $u_i$. Appendix~\ref{appendix:recurrence} shows in
detail how they are uniquely determined by the requirement that $v_i$
and $w_i$ approximate $i^p$ as $i\to\infty$.

Having found the ${u_i}, \, {v_i}$ and $w_i$, the operator $\tilde D$
is given by
\begin{equation}
(\tilde D\Psi)_i=
{8(\tilde\Psi_{i+1}-\tilde\Psi_{i-1})-(\tilde\Psi_{i+2}-\tilde \Psi_{i-2})
\over 12 w_i},
\end{equation}
where we have introduced the shorthand
\begin{eqnarray}
\tilde\Psi_{1/2}&\equiv& v_{1/2}\Psi_{1/2}+u_1\Psi_{3/2}, \\
\tilde\Psi_{3/2}&\equiv& v_{3/2}\Psi_{3/2}+u_1\Psi_{1/2}, \\
\tilde\Psi_i&\equiv& v_i\Psi_i, \qquad i\ge 5/2.
\end{eqnarray}
for the staggered grid and 
\begin{eqnarray}
\tilde\Psi_{0}&\equiv& 0, \\
\tilde\Psi_{1}&\equiv& v_{1}\Psi_{1}+u_{3/2}\Psi_{2}, \\
\tilde\Psi_{2}&\equiv& v_{2}\Psi_{2}+u_{3/2}\Psi_{1}+u_{5/2}\Psi_3, \\
\tilde\Psi_{3}&\equiv& v_{3}\Psi_{3}+u_{5/2}\Psi_2, \\
\tilde\Psi_i&\equiv& v_i\Psi_i, \qquad i\ge 4.
\end{eqnarray}
for the centred grid.

The $\delta_i$ of this method are plotted in Fig.~\ref{deltas_N=2}.  
These $\delta_i$ are uniformly bounded, which demonstrates that this
method is fourth-order accurate.

%%%%%%%%%%%%%%%%%%%%%%%%%%%%%%%%%%%%%%%%%%%%%%%%%%%%%%%%%%%%%%%%%%%%%%%%%

\section{The outer boundary $r=R$}
\label{section:useD0}

We begin by recalling Strand's method \cite{Strand} for treating the
wave equation including boundaries. The one-dimensional wave equation
in first order form is
\begin{equation}
\dot\pi=\psi', \quad \dot\psi=\pi', \quad a\le x\le b
\end{equation}
with energy
\begin{equation}
E=\int_a^b(\pi^2+\psi^2)\,dx, \quad {dE\over dt} = [\pi\psi]_a^b.
\end{equation}
It is natural to discretize this symmetrically in $\pi$ and $\psi$,
that is
\begin{equation}
\dot\Pi = D_0\Psi, \quad \dot\Psi =D_0\Pi,
\end{equation}
with energy
\begin{equation}
\hat E={h\over 2}\left(\Pi^tW_0\Pi+\Psi^tW_0\Psi\right).
\end{equation}
and SBP condition 
\begin{equation}
\label{SBP0}
W_0D_0+(W_0D_0)^t=B_0,
\end{equation}
with with $B_0={\rm diag}(1,0,\dots,0,1)$, as there are two
boundaries. Note that this problem is translation-invariant in the
interior, and so $D_0$ and $W_0$ will naturally be
translation-invariant in the interior, except for finite-sized end
blocks. $D_0$ and $W_0$ with various orders of accuracy in the
interior and at the boundaries have been constructed by Strand
\cite{Strand}. (We have added the suffix $0$ to indicate that this is
the special case $p=0$ of our problem.)

In (\ref{myeqns}) with $p>0$ additional problems result because the
equations are not translation-invariant but depend explicitly on
$r$. In previous Sections we have addressed these problems at interior
points and at the pseudo-boundary $r=0$. 

Strand provides a class of norms $W_0$ that are unit diagonal except
near the boundaries, as well as compatible derivative operators $D_0$
that are the standard minimal width centred difference operators,
except near the boundaries. Hence $D_0$ agrees with our $D$ except at
the outer boundary. Let $W_\infty$ and $\tilde W_\infty$ denote our
previously derived weights for the problem on $0\le r<\infty$, and
simply truncated to the range $i=-M,\dots,M$. We now define operators
with a boundary at $i=\pm M$, corresponding to $r=\pm R$, as follows:
\begin{eqnarray}
W&:=&W_0W_\infty, \\
\tilde W&:=&\tilde W_0W_\infty, \\
D&:=&D_0, \\
\label{Dtildedef}
\tilde D&:=&W^{-1}D_0\tilde W. 
\end{eqnarray}

It is now straightforward to verify that the operators and weights
thus defined obey the desired SBP property (\ref{SBP}) with boundary
operator
\begin{equation}
B=B_0\tilde W_\infty,
\end{equation}
using (\ref{SBP0}). It is essential in this calculation that $\tilde
W_\infty$ and $W_0$ commute. This is true because $\tilde W_\infty$ is
diagonal except near the origin, and $W_0$ is diagonal everywhere and
unit diagonal except near the outer boundary.

The $D$ and $\tilde D$ thus defined agree with their previously
constructed infinitely extended versions except near the boundary, and
so we need to establish their accuracy only near the boundary. By
Strand's construction, using relaxed notation, 
\begin{equation}
h^{-1}D_0={d\over dr}+O(h^{-\tau})
\end{equation}
near the boundary. Also by construction,
\begin{equation}
h^p W_\infty=r^p+O(h^{-2N}), \quad h^p \tilde W_\infty=r^p+O(h^{-2N}),
\end{equation}
near the boundary, with $2N\ge \tau$, and similarly for $\tilde
W_\infty$. Substituting these into (\ref{Dtildedef}), we find
\begin{equation}
\label{Dtildeboundaryaccuracy}
h^{-1}\tilde D={d\over dr}+{p\over r}+O(h^{-\tau}).
\end{equation}
Hence $D$ and $\tilde D$ have the same accuracy both in the interior
and at the boundary, and the same stencil, as the minimal width SBP
operator with diagonal norm $D_0$ of Strand. In this sense, they are
optimal. 

Applying the general prescription above to our second-order accurate
method SBP2 or to the second-order accurate Evans and Sarbach methods
methods reviewed in the Appendix, we have
\begin{eqnarray}
D&=&\left(
\begin{array}{ccccc}
\cdot & & & &  \\
& \cdot & & &  \\
& -{1\over 2} & 0 & {1\over 2} & \\
& & -{1\over 2} & 0 & {1\over 2} \\
& & & -1 & 1 \\ 
\end{array}
\right), \\
\label{Dtildetmp1}
\tilde D &=& \left(
\begin{array}{ccccc}
\cdot & & & & \\
& \cdot & & & \\
& -{v_{M-3}\over 2w_{M-2}} & 0 & {v_{M-1}\over 2w_{M-2}} & \\
& & -{v_{M-2}\over 2w_{M-1}} & 0 & {v_{M}\over 2w_{M-1}} \\
& & & -{v_{M-1}\over w_M} & {v_M\over w_M} \\ 
\end{array}
\right), \\
\label{WWtilde1}
W&=&{\rm diag}\left(\dots, w_{M-2},w_{M-1},{w_M\over
  2}\right), \\
\tilde W&=&{\rm diag}\left(\dots, v_{M-2},v_{M-1},{v_M\over 2}\right), \\
B &=& {\rm diag}\left(\dots,0,0,v_M\right).
\end{eqnarray}
As an example of the general result (\ref{Dtildeboundaryaccuracy}), we
have
\begin{eqnarray}
(\tilde D\Psi)_M&=&-{v_{M-1}\over w_M}\psi_{M-1}+{v_M\over w_M}\psi_M
\nonumber \\
&=&{v_M-v_{M-1}\over w_M}\psi(R)+{v_{M-1}\over w_M}h\psi'(R)
+O(h^2) \nonumber \\
&=&h\left[\psi'(R)+{p\over R}\psi(R)+O(h)\right],
\end{eqnarray}
Hence this method is first-order accurate at the boundary point $i=M$.
The above expressions hold for the Evans, Sarbach and SBP2 methods
with the appropriate $v_i$ and $w_i$. In the last two of these,
$v_i=w_i$.

\begin{widetext}
Applying our general outer boundary prescription to our fourth-order
accurate method SBP4, we can impose accuracy at the boundary of order $\tau
=1$ or $\tau =2$.  For $\tau=1$, following the general prescription
given above, we set
\begin{eqnarray}
D&=&\left(
\begin{array}{cccccc}
\cdot & & & & & \\
& \cdot & & & & \\
& \frac{1}{12} & -\frac{2}{3} & 0 & \frac{2}{3} & -\frac{1}{12} \\
& & \frac{1}{13} & -\frac{8}{13} & 0 & \frac{7}{13} \\
& & & \frac{1}{5} & -\frac{7}{5} & \frac{6}{5} \\ 
\end{array}
\right), \\
\tilde{D}&=&\left(
\begin{array}{cccccc}
\cdot & & & & & \\
& \cdot & & & & \\
& \frac{v_{M-4}}{12w_{M-2}} & -\frac{2v_{M-3}}{3w_{M-2}} & 0 & \frac{2{v}_{M-1}}{3w_{M-2}} & -\frac{{v}_{M}}{12w_{M-2}} \\
& & \frac{v_{M-3}}{13{w}_{M-1}} & -\frac{8v_{M-2}}{13{w}_{M-1}} & 0 & \frac{7{v}_{M}}{13{w}_{M-1}} \\
& & & \frac{v_{M-2}}{5{w}_{M}} & -\frac{7{v}_{M-1}}{5{w}_{M}} & \frac{6{v}_{M}}{5{w}_{M}} \\ 
\end{array}
\right), \nonumber \\ \\
W&=&{\rm diag}\left(\dots, w_{M-2},
   \frac{13{w}_{M-1}}{12},
   \frac{5{w}_M}{12}
\right), \\
\tilde W&=&{\rm diag}\left(\dots, v_{M-2},
   \frac{13{v}_{M-1}}{12},
   \frac{5{v}_M}{12}
\right), \\
B &=& {\rm diag}\left(\dots, 0, 0, v_M\right).
\end{eqnarray}

For $\tau=2$, following the general prescription we set
\begin{eqnarray}
\tilde{D} &=& \left(
\begin{array}{ccccccccc}
\cdot & & & & & & & & \\
& \cdot & & & & & & & \\
& & \frac{v_{M-6}}{12w_{M-4}} & -\frac{2v_{M-5}}{3w_{M-4}} & 0 & \frac{2v_{M-3}}{3w_{M-4}} & -\frac{v_{M-2}}{12w_{M-4}} & 0 & 0 \\
& & 0 & \frac{4v_{M-5}}{49w_{M-3}} & -\frac{32v_{M-4}}{49w_{M-3}} & 0 & \frac{59v_{M-2}}{98w_{M-3}} & 0 & -\frac{3v_{M}}{98w_{M-3}} \\
& & 0 & 0 & \frac{4v_{M-4}}{43w_{M-2}} & -\frac{59v_{M-3}}{86w_{M-2}} & 0 & \frac{59v_{M-1}}{86w_{M-2}} & -\frac{4v_{M}}{43w_{M-2}} \\
& & 0 & 0 & 0 & 0 & -\frac{v_{M-2}}{2w_{M-1}} & 0 & \frac{v_{M}}{2w_{M-1}} \\
& & 0 & 0 & 0 & \frac{3v_{M-3}}{34w_{M}} & \frac{4v_{M-2}}{17w_{M}} & - \frac{59v_{M-1}}{34w_{M}} & \frac{24v_{M}}{17w_{M}}
\end{array}
\right), \\
{W} &=& {\rm diag}\left(\ldots, w_{M-4}, \frac{49w_{M-3}}{48},
 \frac{43w_{M-2}}{48}, \frac{59w_{M-1}}{48}, \frac{17w_{M}}{48}\right) , \\
\tilde{W} &=& {\rm diag}\left(\ldots, v_{M-4}, \frac{49v_{M-3}}{48},
 \frac{43v_{M-2}}{48}, \frac{59v_{M-1}}{48}, \frac{17v_{M}}{48}\right) , \\
B &=& {\rm diag}\left(\dots, 0, 0, v_M\right).
\end{eqnarray}
\end{widetext}
The expression for $D$ is obtained by setting $v_i$ and $w_i$ to $1$
in $\tilde D$. We shall call our SBP4 method with $\tau=1,2$ SBP41 and
SBP42 respectively.

%%%%%%%%%%%%%%%%%%%%%%%%%%%%%%%%%%%%%%%%%%%%%%%%%%%%%%%%%%%%%%%%%%%%%%%%%%

\section{Numerical tests}
\label{section:numerics}

We have implemented our SBP2, SBP41 and SBP42 methods described above,
combined with fourth-order Runge-Kutta (RK4) discretisation in
time. For comparison, we have also implemented the Evans method
(turned into an SBP method by the boundary treatment of
Sec.~\ref{section:useD0}) and the Sarbach method. To complete the
numerical setup, we need to choose continuum boundary conditions at
the physical outer boundary $r=R$ and a way of enforcing them. For our
tests, we choose either homogenous maximally dissipative boundary or
the derivative boundary conditions derived in
Appendixes~\ref{appendix:BCderivslwave} and
\ref{appendix:numericalBC}, and implement them using the Olsson
projection method \cite{Olsson}, which for completeness we review in
Appendix~\ref{appendix:Olsson}.

In SBP4 we use the numerical coefficients $\tilde D_{ij}$, or
equivalently $u_1$, $\bar v_i$ and $\bar w_i$, calculated by the
relaxation method described in Appendix~\ref{appendix:recurrence} up
to $i\sim 2000$, and using the asymptotic results (\ref{asympvbar}),
(\ref{asympwbar}) for larger $i$.

For all evolutions shown here, we use initial data
\begin{equation}
\psi(r,0)=0, \quad \pi(r,0)=e^{-{(r-r_0)^2\over d^2}}+e^{-{(r+r_0)^2\over d^2}}
\end{equation}
with $r_0=5$ and $d=2$. (The Gaussian at negative $r$ is needed to
make $\pi(r)$ strictly even.)  The numerical domain is $0\le r\le R$
with $R=25$. This means that the wave is initially well separated from
both boundaries, and interacts with the symmetry boundary around
$t\sim 5$ and with the outer boundary around $t\sim 20$. We continue
the evolution until $t=40$. 

By construction, all our SBP methods are stable in the energy norm and
consistent with the continuum equations. From the Lax equivalence
theorem we therefore expect convergence to the continuum in the energy
norm $E$, or in other words we expect convergence of $r^{p/2}\pi$ and
$r^{p/2}\psi$ in the uniform $L^2$ norm. We verify this expectation,
but beyond that we also look for pointwise convergence of these
variables.

To check convergence, on a centreed grid we compare evolutions at five
grid resolutions, from $h=1/10$ down to $2^{-4}\cdot 1/10$, each to a
reference evolution at $h=2^{-7}\cdot 1/10$. (By comparison, using
refinement by a factor of 3 on the {\em staggered} grid allows us to
fix $r_M=R$, and while appropriate points of all refined grids still
align with the coarsest grid. Keeping $R$ exactly
resolution-independent is essential for comparing different
resolutions, while aligned grids avoid the need for interpolation.) We
use a Courant factor $\Delta t/h=1/4$ throughout. We plot
\begin{equation}
\label{Epi2def}
e_{\pi,k}(r,t;h)\equiv\left({h\over
1/10}\right)^{-k}r^{p/2}[\pi(r,t;h)-\pi(r,t;h_{\rm ref})]
\end{equation}
and its norm
\begin{equation}
\label{Epi2L2def}
|e_{\pi,k}(\cdot,t;h)|\equiv\left({1\over R}\int_0^R
e_{\pi,k}(r,t;h)^2\,dr\right)^{1/2},
\end{equation}
and similarly for the variable $\psi$.
For $h$ small enough, where a Richardson expansion of the error holds
and is dominated by the leading $O(h^k)$ term, $e_{\pi,k}$ and its norm
should be approximately independent of $h$ (with the differences
generated by subdominant error terms). 

We have tested the Evans, Sarbach, SBP2, SBP41 and SBP42 methods with
a selection of outer boundary conditions and with $p$ in the range
$1\le p\le 22$. Note however the following exceptions: 1. the Evans
method is not defined on the staggered grid for odd $p$; 2. the
Sarbach method is only defined on the centred grid; 3. for $p=1$
Evans, Sarbach and SBP2 on a centred grid are identical. We now
summarise our results. In all evolutions shown in the following {\em
  figures}, we set $p=6$ (corresponding to $l=2$ spherical harmonics
in 3 space dimensions), use a grid centred on $r=0$, and the boundary
condition $\pi=0$.

\begin{figure}
\includegraphics[width=9cm]{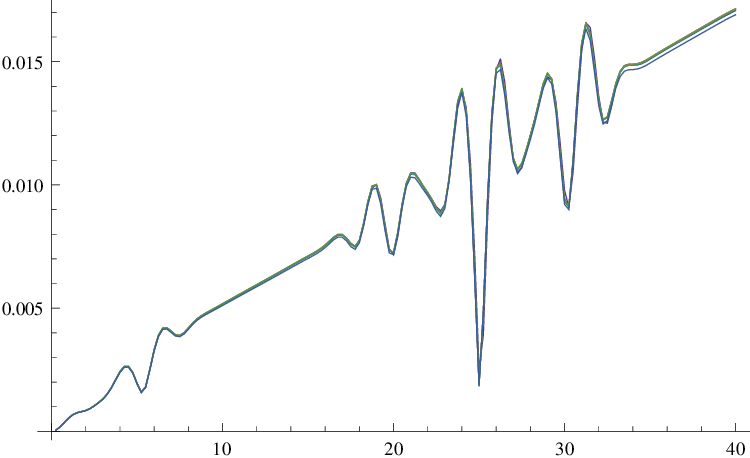}
\caption{\label{figure:SBP2L2} 2nd-order convergence in the energy
  norm of SBP2 for $p=6$, with the initial data given in the text. We
  show $|e_{\pi,2}(\cdot,t;h)|$ against $t$ at 5 different resolutions
  with $h$ decreasing by factors of 2 from $1/10$ to $1/160$. The 5
  curves are on top of each other, demonstrating 2nd-order
  convergence. With the normalisation of Eq.~(\ref{Epi2L2def}), they
  indicate the actual $L^2$ numerical error at resolution $h=1/10$
  (meaning there are $\sim 40$ gridpoints across the wave packet). The
  equivalent curves for $\psi$ and for the Sarbach and Evans numerical
  methods are similar.}
\end{figure}

\begin{figure}
\includegraphics[width=9cm]{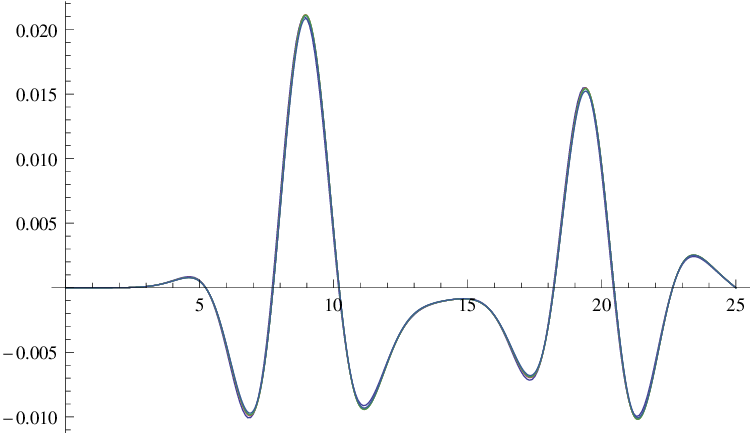}
\caption{\label{figure:SBP2slice} 2nd-order pointwise convergence of
  the same evolution. We show $e_{\pi,2}(r,t;h)$ at $t=14.25$ against
  $r$. The 5 curves are on top of each other, demonstrating perfect
  2nd-order convergence. They indicate the actual pointwise numerical
  error at resolution $h=1/10$. The equivalent curves for $\psi$ and
  for the Sarbach and Evans numerical methods are similar.}
\end{figure}

The three 2nd-order accurate methods, SBP2, Sarbach and Evans all show
2nd-order pointwise convergence (and hence also convergence in the
energy norm) throughout the evolution. Fig.~\ref{figure:SBP2L2}
demonstrates 2nd-order convergence in the energy norm for SBP2, for
all $t$, while Fig.~\ref{figure:SBP2slice} demonstrates pointwise
convergence at $t=14.25$. The error $e_2$ is almost identical for all
methods. Until the wave has interacted with the outer boundary, it
appears smooth, while afterwards there is a small admixture of an
oscillation with the grid frequency.

\begin{figure}
\includegraphics[width=9cm]{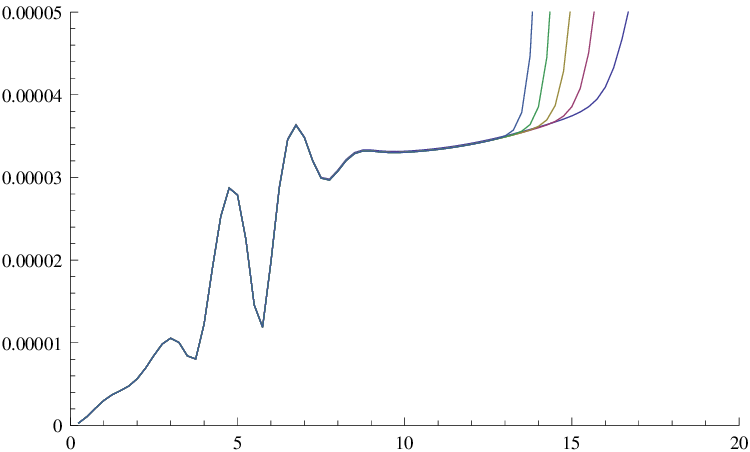}
\caption{\label{figure:SBP41L2k4} 4th-order convergence in the energy
  norm of SBP41. We show $|e_{\pi,4}(\cdot,t;h)|$ against $t$, with all
  other details of the initial data and evolution as for the previous
  figure. The 5 curves are on top of each other, demonstrating perfect
  4th-order convergence until $t\sim 12$, when the interaction of the
  tail of the Gaussian initial data with the outer boundary begins to
  dominate the error. The equivalent curve for $\psi$ looks similar,
  and the equivalent curves for SBP42 are identical until $t\sim 12$.}
\end{figure}

\begin{figure}
\includegraphics[width=9cm]{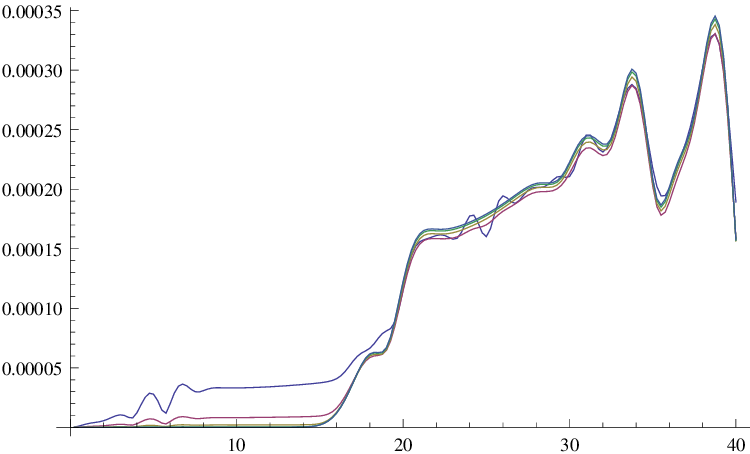}
\caption{\label{figure:SBP41L2k2} 2nd-order convergence in the energy
  norm of SBP41 for $p=6$, {\em after} the wave has first interacted
  with the boundary. We show $|e_{\pi,2}(\cdot,t;h)|$ (instead of
  $e_4$) against $t$, with all other details as in the previous
  figure. The 5 curves are on top of each other, demonstrating
  approximate 2nd-order convergence after $t\sim 20$, when the error
  generated by the interaction of the tail of the Gaussian initial
  data with the outer boundary dominates the error. The equivalent
  curve for $\psi$ looks similar.}
\end{figure}

\begin{figure}
\includegraphics[width=9cm]{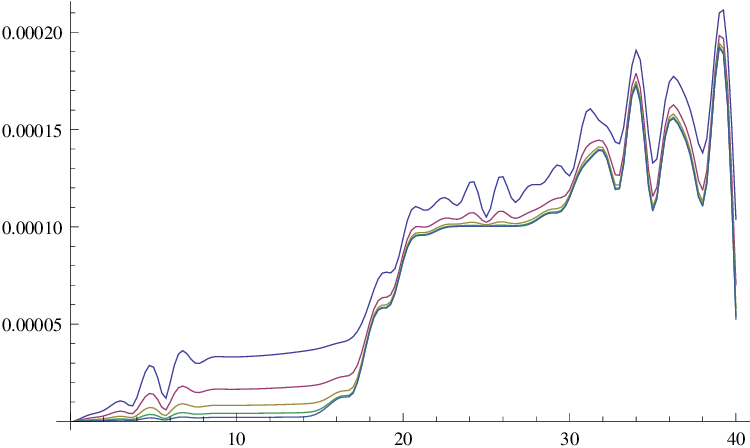}
\caption{\label{figure:SBP42L2k3} 3rd-order convergence in the energy
  norm of SBP42 for $p=6$, after the wave has first interacted with
  the boundary. We show $|e_{\pi,3}(t;h)|$ against $h$, with all other
  details as in the previous figure. The 5 curves are on top of each
  other, demonstrating approximate 3rd-order convergence after $t\sim 20$.}
\end{figure}

\begin{figure}
\includegraphics[width=9cm]{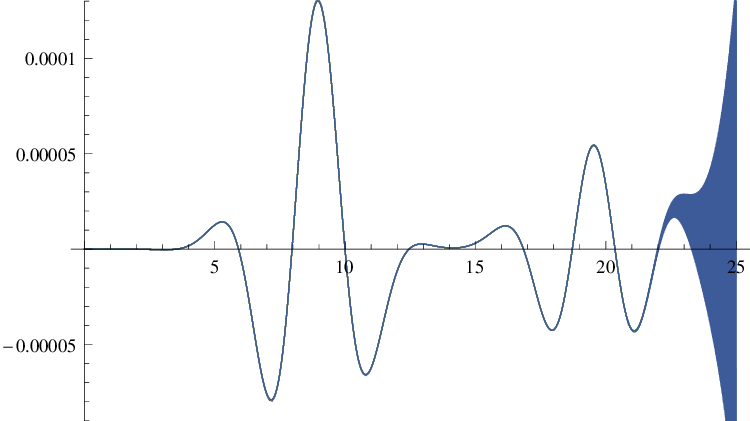}
\caption{\label{figure:SBP42slice} 4th-order pointwise convergence of
  SBP42 before the wave interacts with the outer boundary. We show
  $e_{\pi,4}(r,t;h)$ at $t=14.25$ against $r$ at 5 different
  resolutions. The 5 curves are on top of each other, demonstrating
  pointwise 4-th order convergence, for $0\le r\lesssim 22$, which is
  the region not yet in contact with the outer boundary at this
  time. They again indicate the actual numerical error at resolution
  $h=1/10$. What looks like a filled region for $r\gtrsim 22\le 25$ is
  in fact an oscillation at the grid frequency of the lowest
  resolution, with a smooth envelope. This is the 3rd-order error
  emanating from the outer boundary (and so the curves are no longer
  on top of each other at this resolution). $t=14.25$ has been chosen
  here as the moment when the boundary error is just beginning to
  dominate. Compare Fig.~\ref{figure:SBP2slice}, which shows the same
  moment of time of the evolution with SBP2, with no effect from the
  boundary. The equivalent figures for SBP41 and $\psi$ look similar.}
\end{figure}

For our 4th-order accurate method with 1st and 2nd-order accurate
boundary conditions, SBP41 and SBP42, we see 4th-order pointwise (and
hence energy norm) convergence until the wave interacts with the outer
boundary. Fig.~\ref{figure:SBP41L2k4} demonstrates this for
SBP41. After the wave has interacted with the boundary, SBP41 drops to
2nd-order convergence in the energy norm (see
Fig.~\ref{figure:SBP41L2k2}), while SBP42 drops to 3rd-order
convergence in the energy norm (see Fig.~\ref{figure:SBP42L2k3}). Note
that in each case the global accuracy is one order higher than the
accuracy $\tau$ of $D$ and $\tilde D$ at the boundary. For both
methods, the error after the interaction with the boundary is
dominated by an oscillation with the grid frequency, with a smooth
envelope, and so they do not converge pointwise in the standard sense,
although the envelope of the grid frequency noise
does. Fig.~\ref{figure:SBP42slice} is a snapshot that shows the
transition from 4th-order pointwise convergence to this behaviour as
the wave begins to interact with the boundary.

Two comments on our convergence tests are worth making: First, note
that the Lax theorem only gives convergence in $L^2$ of $r^{p/2}\pi$,
$r^{p/2}\psi$. We do find this in our tests, but we also find
pointwise convergence at the same rates, at all times for SBP2, and
for SBP4 before the wave interacts with the outer boundary. One can go
further and look at the convergence of the unscaled variables $\pi$,
$\psi$, for which the theory makes no prediction. We find that they
converge pointwise for $p\lesssim 4$ at all times, and for all $p$
while the wave is away from the symmetry boundary. However, while the
wave is reflected at the origin, the continuum solution oscillates
rapidly approximately $p$ times. (This can be shown by constructing
the exact solution as a sum involving the first $p$ derivatives of the
initial data.) A small phase error at this stage gives rise to a very
large pointwise error and pointwise convergence is lost (at the
resolutions we ran). However, as the solution moves out again,
different resolutions agree again much better. This is
compatible with the observed pointwise convergence of $r^{p/2}\pi$,
$r^{p/2}\psi$ because in these rescaled variables the complicated
continuum behaviour at the origin is hidden and so is the momentary
increase of the error.

\begin{figure}
\includegraphics[width=9cm]{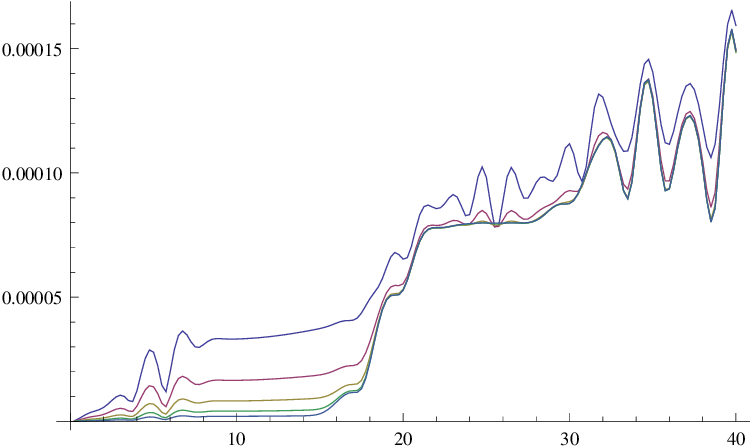}
\caption{\label{figure:SBP42L2k3_derivBC} 3rd-order convergence in the
  energy norm of SBP42 for $p=6$, with a {\em derivative} boundary
  condition, after the wave has first interacted with the
  boundary. This is similar to Fig.~\ref{figure:SBP42L2k3}, except
  that the boundary condition at $r=R$ is now $\pi+\pi'=0$ instead of
  $\pi=0$.}
\end{figure}

Our second comment is that the general theory for the accuracy of
first-order hyperbolic initial-boundary value problems
\cite{GKO,SvardNordstrom,Gustafsson1975,Gustafsson1981} suggests that
the order of global accuracy is determined by the lower of the order
of the physical boundary conditions, and the order of purely numerical
(``extra'') boundary conditions plus 1. In our case we always have one
physical and one extra boundary condition. The accuracy order of the
physical boundary condition is $\infty$ for maximally dissipative
physical boundary conditions ($\mu=\nu=0$) and $\tau$ for physical
boundary conditions involving a derivative ($\mu$ or $\nu\ne 0$), as
we discretise these using $D$ and $\tilde D$ on the boundary. The
accuracy order of the extra boundary condition is always $\tau$, as it
relies on evaluating $D$ and $\tilde D$ on the boundary. Hence we
would expect global accuracy of order $\tau+1$ for any maximally
dissipative boundary condition, and $\tau$ for any boundary condition
involving a derivative. However, experimentally we find $\tau+1$ in
both cases, which means that the maximally dissipative boundary
conditions perform as expected, and the derivative boundary conditions
perform one order better than expected. The latter point is illustrated in
Fig.~\ref{figure:SBP42L2k3_derivBC}.

%%%%%%%%%%%%%%%%%%%%%%%%%%%%%%%%%%%%%%%%%%%%%%%%%%%%%%%%%%%%%
%%%%%%%%%%%%%%%%%%%%%%%%%%%%%%%%%%%%%%%%%%%%%%%%%%%%%%%%%%%%%

%%%%%%%%%%%%%%%%%%%%%%%%%%%%%%%%%%%%%%%%%%%%%%%%%%%%%%%%%%%%%
%%%%%%%%%%%%%%%%%%%%%%%%%%%%%%%%%%%%%%%%%%%%%%%%%%%%%%%%%%%%%

%%%%%%%%%%%%%%%%%%%%%%%%%%%%%%%%%%%%%%%%%%%%%%%%%%%%%%%%%%%%%
%%%%%%%%%%%%%%%%%%%%%%%%%%%%%%%%%%%%%%%%%%%%%%%%%%%%%%%%%%%%%

%%%%%%%%%%%%%%%%%%%%%%%%%%%%%%%%%%%%%%%%%%%%%%%%%%%%%%%%%%%%%
%%%%%%%%%%%%%%%%%%%%%%%%%%%%%%%%%%%%%%%%%%%%%%%%%%%%%%%%%%%%%

%%%%%%%%%%%%%%%%%%%%%%%%%%%%%%%%%%%%%%%%%%%%%%%%%%%%%%%%%%%%%
%%%%%%%%%%%%%%%%%%%%%%%%%%%%%%%%%%%%%%%%%%%%%%%%%%%%%%%%%%%%%

%%%%%%%%%%%%%%%%%%%%%%%%%%%%%%%%%%%%%%%%%%%%%%%%%%%%%%%%%%%%%
%%%%%%%%%%%%%%%%%%%%%%%%%%%%%%%%%%%%%%%%%%%%%%%%%%%%%%%%%%%%%

%%%%%%%%%%%%%%%%%%%%%%%%%%%%%%%%%%%%%%%%%%%%%%%%%%%%%%%%%%%%%
%%%%%%%%%%%%%%%%%%%%%%%%%%%%%%%%%%%%%%%%%%%%%%%%%%%%%%%%%%%%%

As a further test of the predicted behaviour of our methods, we have
also evaluated the discrete energy at every time step. With the
boundary conditions $\rho\pi+\mu\pi'=0$ and
$\sigma\psi+\nu(\psi'+p\psi/r)=0$ discussed in
Appendix~\ref{appendix:numericalBC} we have $d\hat E_b/dt=0$. With the
maximally dissipative boundary condition $\rho\pi+\sigma\psi=0$ we
have $d\hat E/dt=\chi\pi_M\psi_M\le 0$, and we have evolved the
expected value of $\hat E$ by discretising this in $t$ using RK4.  In
all these cases the discrepancy between the evaluated and predicted
numerical energy is of relative size $10^{-8}$, essentially
independent of the choice of SBP method and of the resolution, and
increases linearly with $t$. These observations are compatible with
the expectation of accumulated round-off (machine precision) error,
with zero finite-differencing error.

The energy of SBP4 is not positive definite on the staggered grid for
$p=1,2$, and so we would not expect it to be stable. However, we do
not see signs of instability in our numerical experiments.

The Sarbach method behaves like SBP2 and Evans for $p\lesssim 8$, but
requires a much smaller Courant factor in order to be stable for
larger $p$: for $p=10$, $12$ and $22$, we empirically find that the
Courant factor needs to be reduced to $1/8$, $1/16$ and $1/800$,
respectively. By contrast, all other SBP methods are stable with RK4
with a Courant factor of $1/4$ up to $p=22$.

Finally, we have also implemented the obvious naive, non-SBP,
second-order accurate finite difference method in which all
derivatives are just evaluated using centred derivatives, and the
$p/r$ term is evaluated pointwise, assuming a staggered grid. In our
notation this corresponds to defining
\begin{equation}
\label{obvious}
\tilde D\psi_i=D\psi_i+p\psi_i/r_i
\end{equation}
on a staggered grid, with $D$ given by (\ref{D2}) and ghost points at
the origin. This method is unstable at the origin for all $p>0$, with
blowup occurring more rapidly for larger $p$, and more rapidly at
higher resolution. This failure of the ``standard'' method (which is
SBP and hence stable for $p=0$) was of course the motivation for our
work. [At the outer boundary $r=R$, we implemented ``copy''
  (zeroth-order extrapolation) boundary conditions for this test, but
  we moved the outer boundary very far out so that the wave does not
  interact with the boundary, even numerically, before the blowup
  occurs. We are therefore certain that the instability of this method
  is due to the $p/r$ term and not to our particular choice of outer
  boundary condition.] 

%%%%%%%%%%%%%%%%%%%%%%%%%%%%%%%%%%%%%%%%%%%%%%%%%%%%%%%%%%%%%%%%%%%%%%%%%%

\section{Conclusions}
\label{section:conclusions}

It is surprising that the lower-order term $p/r$ in (\ref{myeqns})
alone can make standard finite differencing schemes unstable, and that
an elaborate SBP scheme is necessary. Note however that a standard
centred finite difference implementation of the one-dimensional wave
equation is already SBP except possibly at the boundaries, while the
equivalent naive finite differencing of (\ref{myeqns}) for $p>0$ is
not SBP even at interior points.

It seems highly unlikely to us that any scheme for (\ref{myeqns}) that
is not SBP can be made stable without using numerical dissipation, for
any choice of discrete boundary condition. Numerical dissipation {\em
  can} in fact stabilise the non-SBP discretisation (\ref{obvious}),
but more and more dissipation is required with increasing $p$, making
this approach useless for even moderately large $p$. Again we suspect
that this will be so for any non-SBP scheme. This rules out non-SBP
finite differencing schemes for large $p$. Furthermore, in
applications where the physical growth or decay of the continuum solution
is under investigation (for example, in stellar perturbation
theory), the numerical method should be as little dissipative as
possible.

The Evans method has been used with success previously (see
\cite{HondaChoptuik} for a $p=2$ application and \cite{critpert} for
$p\ge 2$), but we have here turned it into a complete SBP method by
the appropriate modification at the outer boundary $r=R$. This
modification would not have been obvious outside of an SBP
framework. The Evans method and our SBP2 method work equally well for
all $p$. The Evans method is simpler to implement, but it does not
exist for odd $p$ on a centred grid, in which case SBP2 can be used
instead. 

For higher accuracy, our SBP42 method should be used. It works for any
$p$ on both centred and staggered grid. (SBP41 is described here only
for presentation purposes and numerical tests). It requires loading
the coefficients of $\tilde D$ from a file \cite{data}, but is
otherwise as simple to implement as any other method, and its stencil
has only 5 points (except near the origin), as narrow as possible for
a 4th-order accurate method.

In hindsight we note that discrete energy conservation bounds only
$r^{p/2}\pi$ and $r^{p/2}\psi$. As the maximum of $r^{-p/2}$ on the
grid increases as $h^{-p/2}$ with resolution, the maximum of the
numerical solution can in principle increase by the same factor,
allowing it to become very much larger than the continuum solution. We
find empirically that this happens in the Sarbach method with
$p\gtrsim 8$, effectively leading to blowup even though a numerical
energy $\hat E$ is conserved, unless the Courant number is severely
reduced, but that it does not happen in the other SBP methods. We have
no rigorous explanation for this, but it may be connected to the fact
that the local error near the origin in the Sarbach method is
dominated by $O(h^2/r^2)=O(i^2)$ terms while the local error in the
other methods is $O(h^2)$ uniformly in $r$.

The construction of our SBP4 method is designed to achieve a {\em
  uniform} in $r$ bound on the local error (and the failure of the
Sarbach method at large $p$ seems to justify the need for this). It
may be possible that a uniformly fourth-order accurate SBP method
exists in which the coefficients of $\tilde D$ can be given in closed
form (as they are for the second-order accurate Evans method), but we
have not found such a method.

To summarise our results: Until now, the only known stable numerical
method for the wave equation (\ref{myeqns}) on the semi-infinite
domain $0\le r<\infty$ was the Evans method. We have shown that it is
stable because it is SBP. We have generalised it to the finite domain
$0\le r\le R$, and to grids both centred and staggered with respect
to $r=0$, for arbitrary $p$.

Going beyond 2nd-order accuracy, we have given 4th-order accurate SBP
operators on this finite domain on both centred and staggered grids,
and we have described a general strategy for constructing SBP
operators of {\em arbitrary} accuracy. We have proved SBP for these
methods for the usual maximally dissipative boundary conditions at
$r=R$, which include Dirichlet and Neumann boundary conditions, and
for two families of boundary conditions involving first derivatives of
$\pi$ or $\psi$. 

Our work can be seen as generalising the work of Strand on SBP
operators of arbitrary accuracy from the case $p=0$ to the case $p>0$,
motivated by applications of the wave equation in spherical rather
than Cartesian coordinates.

%%%%%%%%%%%%%%%%%%%%%%%%%%%%%%%%%%%%%%%%%%%%%%%%%%%%%%%%%%%%%%%%%%%%%%%%%%

\acknowledgments

We would like to thank Olivier Sarbach for instructive conversations
and comments on the manuscript, and Piotr Biz\'on
for pointing out the application of our methods to general $n$.
JMM was supported by ANR grant BLAN07-1\_201699 
``LISA Science'', and also in part by MICINN projects
FIS2009-11893 and FIS2008-06078-C03-03. CG would like to thank
GReCo/IAP and LUTH/Observatoire de Meudon for hospitality, and was
partly supported by ANR grant 06-2-134423 ``Mathematical Methods in
General Relativity''. DG was supported by NSF
grant PHY-0855532.

%%%%%%%%%%%%%%%%%%%%%%%%%%%%%%%%%%%%%%%%%%%%%%%%%%%%%%%%%%%%%%%%%%%%%%%%%%

\appendix

%%%%%%%%%%%%%%%%%%%%%%%%%%%%%%%%%%%%%%%%%%%%%%%%%%%%%%%%%%%%%%%%%%%%%%%%%%%%

\section{Rigorous treatment of ghost points at $r=0$}
\label{appendix:ghostpoints}

We initially assume a staggered grid.
Consider $\dot \Psi_i$ for physical grid points
$i>0$. We can write the use of ghost points explicitly as
\begin{equation}
\dot \Psi_i=h^{-1}\sum_{j>0} \left(D_{ij}\Pi_j+D_{i,-j}\Pi_{-j}\right)
=h^{-1}\sum_{j>0}D^{(+)}_{ij}\Pi_j, 
\end{equation}
where
\begin{equation}
\label{D+}
D^{(+)}_{ij}\equiv D_{ij}+D_{i,-j}, \qquad i,j>0.
\end{equation}
We think of this as ``folding over the ghost points''. A similar
observation holds for $\tilde D$, except that as $\Psi_i$ is odd, the
equivalent of (\ref{D+}) is
\begin{equation}
\label{Dtilde+}
\tilde D^{(+)}_{ij}\equiv \tilde D_{ij}-\tilde D_{i,-j}, \qquad i,j>0.
\end{equation}
Note that $D^{(+)}_{ij}\ne D_{ij}$ even for $i,j>0$, thus requiring a
separate symbol. (The symbol $D^{(+)}$ is a reminder of the range
$i,j>0$.) The split of $D^{(+)}_{ij}$ into $D_{ij}$ and $D_{i,-j}$ for
$i,j>0$ is in general not unique. We do, however, have a natural
prescription for this split if we assume that $D_{ij}$ is
translation-invariant, i.e. depends only on $i-j$ even at the
boundary.

In order to extend $D_{ij}$ and $\tilde D_{ij}$ to negative $i$, we
use the requirement that (\ref{discreteevenodd}) hold at all times, or
\begin{equation}
\label{discreteevenodddot}
\dot\Pi_{-i}=\dot\Pi_i, \qquad \dot\Psi_{-i}=-\dot\Psi_i.
\end{equation}
The first equation of (\ref{discreteevenodd}) and the second equation
of (\ref{discreteevenodddot}) immediately give, for $i,j>0$, that
\begin{equation}
\label{tmp1}
D_{ij}+D_{i,-j}+D_{-ij}+D_{-i,-j} = 0 .
\end{equation}
The second equation of (\ref{discreteevenodd}) and the first equation
of (\ref{discreteevenodddot}), after substituting (\ref{SBPbis}) and
using (\ref{Wsymmetry}), give, for $i,j>0$, that
\begin{equation}
(W^{-1})_{ik}\left(D_{kl}-D_{k,-l}-D_{-kl}+D_{-k,-l}\right)W_{lj} = 0,
\end{equation}
and hence
\begin{equation}
\label{tmp2}
D_{ij}-D_{i,-j}-D_{-ij}+D_{-i,-j} = 0 .
\end{equation}
Taking the sum and difference of (\ref{tmp1}) and (\ref{tmp2}), we
obtain (\ref{Dsymmetry}).

Finally, we have, for $i,j>0$, that
\begin{eqnarray}
\tilde D^{(+)}_{ij} &=& \sum_{k,l>0}
-(W^{-1})_{ik}(D_{lk}-D_{-lk})\tilde W_{lj} 
\nonumber \\
&=& \sum_{k,l>0}
-(W^{-1})_{ik}(D_{lk}+D_{l,-k})\tilde W_{lj} \nonumber \\
&=& \sum_{k,l>0}
-(W^{-1})_{ik}(D^{(+)t})_{kl}\tilde W_{lj},
\end{eqnarray}
and so (\ref{SBP}) holds for the operators $D^{(+)}$ and $\tilde D^{(+)}$ with
ghost points folded in if and only if it holds for the extened
operators $D$ and $\tilde D$. This confirms that the
introduction of ghost points is just a matter of notation (or coding).

A similar argument goes through on a centred grid, with the point
$i=0$ ``split'' between the domains $r\ge 0$ and $r\le 0$. Here we
note only that when removing the ghost points, the discrete energy on
a centred grid is
\begin{eqnarray}
\hat E&=&{h^{p+1}\over 4}(W_{00}\Pi_0^2+\tilde W_{00}\Psi_0^2)
\nonumber \\ &&
+{h^{p+1}\over 2}\sum_{i,j=1}^M (W_{ij}\Pi_i\Pi_j+\tilde W_{ij}\Psi_i\Psi_j).
\end{eqnarray}
(Note the $1/4$.)

%%%%%%%%%%%%%%%%%%%%%%%%%%%%%%%%%%%%%%%%%%%%%%%%%%%%%%%%%%%%%%%%%%%%%%%%%%

\section{Solution of the recurrence relations}
\label{appendix:recurrence}

In the $N=1$ case, the recurrence relation
Eq.~(\ref{SBP2recurrence}) for the $w_i$ is solved as follows. To
work with a bounded quantity, we define the new variable $\bar{w}_i$
as
\begin{equation}
\label{wbar}
w_i = i^p \bar{w}_i .
\end{equation}
(Therefore $\bar w_{-i}=-\bar w_i$ for odd $p$, while $w_i>0$.)
It obeys the linear recurrence relation
\begin{equation}
\bar{w}_i =
 \frac{2(p+1)}{i}\left(1-\frac{1}{i}\right)^p \bar{w}_{i-1}
 + \left(1-\frac{2}{i}\right)^{p+1} \!\!\!\! \bar{w}_{i-2} .
\end{equation}
Trying asymptotic solutions of the form 
\begin{equation}
\label{asymptotic}
\bar{w}_i = \rho^i \sum_{m=0}^\infty C_m i^{k-m}
\end{equation}
for constants $\rho$ and $k$ shows that the two linearly independent
solutions have $\rho=\pm 1$ and are (fixing a constant overall factor,
and assuming $i>0$)
\begin{eqnarray}
\bar{w}^{(+)}_i &=& 1 + \frac{p(p^2-1)}{12i^2}+O (i^{-4}) ,
\label{SBP2mode+} \\
\bar{w}^{(-)}_i &=& (-1)^i i^{-2(p+1)}
 \Biggl( 1 -\frac{(p+1)(p+2)(p+3)}{12i^2}\nonumber \\   && 
\qquad\qquad\qquad\qquad+ O (i^{-4}) \Biggr) .
\label{SBP2mode-}
\end{eqnarray}
The first one is asymptotically constant, and the second is an
oscillating decaying solution. The general asymptotic solution is an
arbitrary linear combination of those, and hence it is also
asymptotically constant.  
[The asymptotically constant mode (\ref{SBP2mode+}) for given
$p$ is a finite polynomial in $i^{-2}$ of (the integer part of)
$1+p/2$ terms.  For example, restricting to $i>0$, for $p=1$ we have
$\bar{w}_i^{(+)} = 1$ and for $p=2$ we have
$\bar{w}_i^{(+)} = 1 + 1/(2i^2)$.]

From these results we can infer the asymptotic behaviour of the
$\delta_{\alpha i}$.  If $\bar{w}_i$ tends to a constant
$\bar{w}_\infty\not= 0$ then we have
\begin{eqnarray}
\delta_{0i} &=& \frac{p(1-p)}{2} + O(i^{-2}) , 
\label{asymp1} \\
\delta_{1i} &=& \frac{p}{2} + O(i^{-2}) .
\label{asymp2}
\end{eqnarray}
Only in the case where $\bar{w}_\infty = 0$ and only the oscillating
mode is present is there a divergence in $\delta_{1}$, namely
\begin{eqnarray}
\delta_{0i} &=& 2i^2 + \frac{(p+2)(p+3)}{2} + O(i^{-2}) , \\
\delta_{1i} &=& \frac{p+2}{2} + O(i^{-2}).
\end{eqnarray}
However, with our initial data $w_{-1/2}=w_{1/2}$
on the staggered grid or $w_1=(1+p)w_0$ on the centred grid, the
constant solution is present, and hence the
sequence $\delta_{\alpha i}$ converges as $i\to\infty$, and is
therefore bounded. Furthermore, the upper bound of its absolute value
is close to the asymptotic value, as we show in Fig.~\ref{deltas_N=1}.

Finally, we adjust the arbitrary overall factor such that $\lim_{i\to
\infty}\bar w_i=1$. On the centred grid we need
\begin{equation}
w_0 = \frac{p!}{2^p}
\end{equation}
for any value of $p$. $\bar w_i$ for $i>0$ is then actually
given by the asymptotically constant polynomial
(\ref{SBP2mode+}). (For even $p$, this is true for all $i$, but not for
$i=0$ with odd $p$, where the special form of the accuracy condition
at the centre needs to be used.)

On the staggered grid we need for even $p$
\begin{equation}
\bar{w}_{1/2} = \frac{[(p+1)!!]^2}{p+1} ,
\end{equation}
which also leads to the polynomials (\ref{SBP2mode+}). However for odd
$p$ the symmetry condition at the centre is incompatible with having
only the asymptotically constant mode, and we need a contribution from
the oscillating mode (\ref{SBP2mode-}).  For $\lim_{i\to \infty}\bar
w_i=1$ we now need
\begin{equation}
\bar{w}_{1/2} = \frac{2}{\pi}\frac{[(p+1)!!]^2}{p+1}.
\end{equation}

The method for solving the recurrence relation in the $N=2$ case is similar.
With the equivalent of (\ref{wbar}) and (\ref{asymptotic}) for $v_i$,
the fourth-order linear recurrence for $\bar{v}_i$ has four
independent asymptotic solutions with
\begin{equation}
\rho=1,\quad -1,\quad 4+\sqrt{15},\quad 4-\sqrt{15},
\end{equation}
all with $k=0$.  The linearity of the recurrence relation implies that
the general solution $\bar{v}_i$ is a linear combination of the four
corresponding modes $\bar{v}^{(\rho)}_i$. It is possible to show that
if the linear combination contains any contribution of the growing or
oscillating modes then the $\delta_\alpha$ are not bounded. Hence we
must find a solution which only contains the asymptotically constant
and the decaying modes.  The freedom in $u_1$, $v_{1/2}$ and $v_{3/2}$
on the staggered grid, and in $u_{3/2}$, $u_{5/2}$ and $v_1$ on the
centred grid allows us precisely to cancel simultaneously the growing
mode and the oscillating non-decaying mode and fix an overall constant
factor. To do that we proceed as follows.

We first compute three arbitrary solutions of the recurrence up to
some high value of $i$, say 1000. For example, on the staggered grid we
can set each of $u_1$, $v_{1/2}$ and $v_{3/2}$ to 1 and the other two to 0.
The three solutions are dominated by the growing mode, and reach
very high values, of order $(4+\sqrt{15})^{1000}\sim 10^{896}$.
We have detected extreme sensitivity of the solution to the initial
conditions, roughly losing one decimal digit of precision per iteration,
and hence the recurrence is solved with exact rational arithmetic,
using {\em Mathematica}.

Then we compute the asymptotic form of the modes, up to order $O(i^{-8})$.
For instance for the asymptotically constant mode we have
\begin{eqnarray}
\bar{v}_i^{(1)} = 1
&+& \frac{(2p-1)(p-1)p(p+1)(p+3)}{60i^4} \nonumber \\
&+& \frac{(2p-3)(p-3)(p-2)(p-1)p(p+1)(p+3)}{504i^6} \nonumber \\
&+& O(i^{-8}).
\label{asympvbar}
\end{eqnarray}
[for $p=2$ this is simply $1+3/(2i^4)+O(i^{-8})$.]
However, in contrast to the $N=1$ case, these are finite
polynomials only for odd $p$, but not for even $p$.
For $i\sim 1000$ this expression will give results correct up to relative
errors smaller than $10^{-18}$ for $p\le 10$.  We take three such
values of $i$ and construct a linear system to find which linear
combination of our three solutions gives that mode
$\bar{v}_i^{(1)}$. For such high values of $i$ we can neglect the
contribution of the decaying mode. In this way we determine the values
of $\bar{v}_i$ up to $i=1000$.  For larger $i$, and $p\le 10$, the
asymptotic series are accurate to 16 digits.  In our experiments below
we shall use up to $p=22$, for wich values up to $i=2000$ must be
computed to use the given asymptotic expansions with relative errors
below double precision.  Note that we do not know if these series are
convergent.

From $v_i$ we can compute $w_i$. 
This gives the following asymptotic behaviour for $\bar{w}_i$,
\begin{eqnarray}
\label{asympwbar}
\bar{w}_i = 1
&+& \frac{(2p+1)(p+1)p(p-1)(p-3)}{60i^4} \nonumber \\
&+& \frac{(2p-1)(p-5)(p-3)(p-2)(p-1)p(p+1)}{504i^6} \nonumber \\
&+& O(i^{-8}) .
\end{eqnarray}

For $i\ge 9/2$ on the staggered grid and $i\ge 6$ on the centred grid
the $\delta$ can be computed from $v_i$ and $w_i$ as follows,
\begin{eqnarray}
\delta_{0i} &=& i^5\,\frac{v_{i-2}-8v_{i-1}+8v_{i+1}-v_{i+2}}{12w_i} -p \,i^4 , \label{delta1fromvw} \\
\delta_{1i} &=& i^2\,\frac{v_{i-2}-v_{i-1}-v_{i+1}+v_{i+2}}{9w_i} , 
\label{delta2fromvw} \\
\delta_{2i} &=& i \,\frac{2v_{i-2}-v_{i-1}+v_{i+1}-2v_{i+2}}{36w_i} .
\label{delta3fromvw}
\end{eqnarray}
The previous expansions imply that the $\delta$ are bounded and have
finite limits:
\begin{eqnarray}
\delta_{0i} &=& -p(p-1)^2 + O(i^{-2}) \label{asymp1b} \\
\delta_{1i} &=& \frac{p(p-1)}{3} + O (i^{-2}), \label{asymp2b} \\
\delta_{2i} &=& -\frac{p}{6} + O (i^{-2}) . \label{asymp3b}
\end{eqnarray}
The limit value of $\delta_0$ is cubic in $p$. That means that
$\delta_0$ is very large for large values of $p$. Comparing with $N=1$
it is plausible that $\delta_0$ has an asymptotic limit which grows
like $p^{N+1}$.

We provide in our webpage \cite{data} data files with double-precision
results for $\bar{v}$, $\bar{w}$ and the $\delta$ for $1\le p\le 22$
and $i\le 2000$. Formulas (\ref{asympvbar}) and (\ref{asympwbar}) can
be used to compute $\bar{v}$ and $\bar{w}$ for these $p$ and $i>2000$
to 16 digits.

We find that for $p=1$ and $p=2$ (the wave equation in cylindrical and
spherical symmetry), on the staggered grid, $v_{1/2}<0$, so that
$\tilde W$ is then not positive definite. This problem is absent for
$p\ge 3$, or on the centred grid. It is possible that allowing for
$u_i$ other than $u_1$ to be nonzero this could be fixed, but we have
not tried this. 

%%%%%%%%%%%%%%%%%%%%%%%%%%%%%%%%%%%%%%%%%%%%%%%%%%%%%%%%%%%%%%%%%%%%%%%%%%

\section{The Evans method}
\label{appendix:Evans}

\begin{figure*}[ht!]
\psfrag{ipoint}{\large $i$}
\psfrag{d1}{\large $\!\!\!\delta_0$}
\psfrag{d2}{\large $\!\!\!\delta_1$}
\includegraphics[width=18cm]{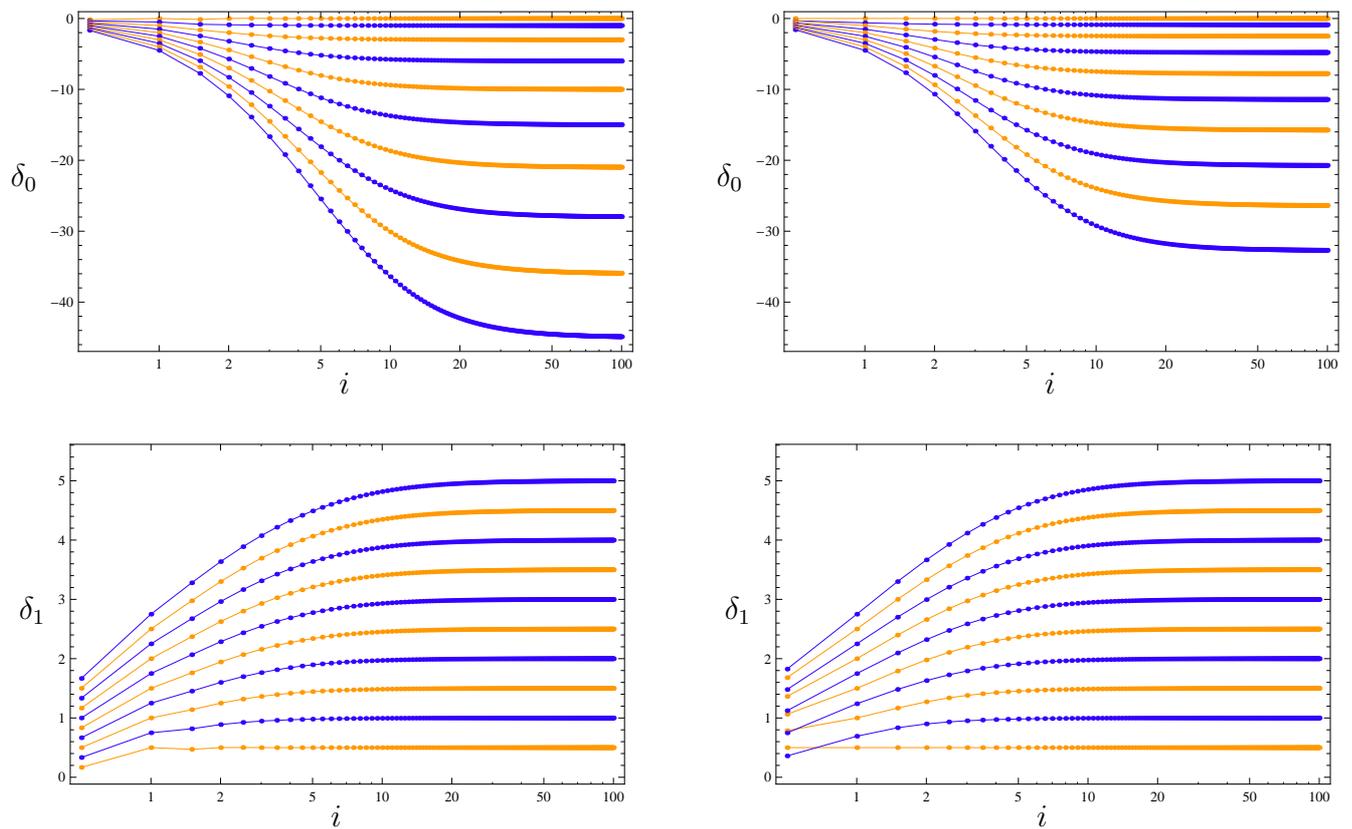}
\caption{\label{deltas_N=1} Values of $\delta_0$ and $\delta_1$ for
  our second-order accurate ($N=1$) methods, for $p=1,\dots,10$. SBP2
  is in the left column and Evans in the right column. The staggered
  grids (half-integer $i$) and centred grids (integer $i$) are shown
  on the same plot. In all cases increasing values of $p$ correspond
  to increasing $|\delta_i|$, with even values of $p$ shown in blue
  (dark) and odd values in orange (light). Note that the Evans method
  does not exist for odd $p$ on the centred grid, and the
  corresponding dots are absent. We see a rapid convergence towards
  the respective asymptotic values (\ref{asymp1}) and (\ref{asymp2})
  for SBP2, and $\delta_0 \to (p+2)p(1-p)/3/(p+1)$ and $\delta_1\to
  p/2$ for the Evans method.  Note that the $i$ axis is logarithmic.
}
\end{figure*}

Here we review the method of Evans \cite{Evans} in the notation of our
paper and present a boundary treatment that makes it SBP. The
continuum identity
\begin{equation}
\psi'+\frac{p}{r}\psi =
(p+1)\frac{d(r^p\psi)}{d(r^{p+1})}
\label{Evansidentity}
\end{equation}
suggests the difference operator $\tilde D$ given by
\begin{equation}
\label{Evans1}
h^{-1}(\tilde
D\Psi)_i=(p+1){r_{i+1}^p\Psi_{i+1}-r_{i-1}^p\Psi_{i-1}\over
  r_{i+1}^{p+1}-r_{i-1}^{p+1}}.
\end{equation}
We combine it with the usual second-order accurate 3-point symmetric
difference operator (\ref{D2}). Comparing (\ref{Dtilde2}) with
(\ref{Evans1}) we see that the Evans method is then SBP with
\begin{eqnarray}
\label{Evansv}
v_i &=&i^p, \\
\label{Evansw}
w_i& =&{(i+1)^{p+1}-(i-1)^{p+1} \over 2(p+1)}.
\end{eqnarray}
We note that these are well defined for all $i$ including
$i=0$. Indeed, the accuracy conditions at the origin
(\ref{cequalities_origin}) for a method with diagonal energy
(\ref{vwdef}) reduce to $v_1=(1+p)w_0$, which is easily seen to hold
for the ansatz (\ref{Evansv},\ref{Evansw}) for even $p$. Note that the
Evans method does not work for odd $p$ (the wave equation in even
space dimensions) on a staggered grid, as then $w_0=0$. The plots in
the right half of Fig.~(\ref{deltas_N=1}) show that this method is
second-order accurate uniformly in $r$, like our method SBP2. 

To our knowledge, no SBP treatment of the outer boundary for the SBP
method has been given. However, our general method of
Sec.~\ref{section:useD0} immediately gives us a prescription, namely 
(\ref{WWtilde1},\ref{Dtildetmp1}). 

The identity (\ref{Evansidentity}) seems at first sight to 
suggest a generalization of the
Evans method to accuracy order $2N$, discretizing $d/dr+p/r$ as
\begin{equation}
\label{Evansgeneralized}
(p+1)\frac{D(r^p\Psi)}{D(r^{p+1})},
\end{equation}
where $D$ is some discretization of $d/dr$ of accuracy order $2N$.
If we take the norms
\begin{equation}
v_i = i^p,
\qquad
w_i = \frac{1}{p+1}D(i^{p+1}),
\end{equation}
then (\ref{Evansgeneralized}) also obeys the SBP property. However,
for the minimal-width centred stencils $D$ of order larger than 2
this does not work. To see this we differentiate $\psi(r)=a r+b
r^3$ using the operators $D$ of accuracy $N=1$ and $N=2$.
The discretization errors are, respectively,
\begin{eqnarray}
 b (p+3)(2p+1) h^2 &+& O(h^3), \\
- \frac{2b (p+3)p(p-1)(2p-1)}{15} \frac{h^4}{r^2} &+& O(h^5).
\end{eqnarray}
In the latter case we see an error of the form $h^4/r^2$, which becomes
$h^2$ near the centre.

We have not been able to generalize the Evans method to avoid this type
of singular error term.

%%%%%%%%%%%%%%%%%%%%%%%%%%%%%%%%%%%%%%%%%%%%%%%%%%%%%%%%%%%%%%%%%%%%%%%%%%

\section{The Sarbach method}
\label{appendix:Sarbach}

Here we review the method of Sarbach \cite{CalabreseNeilsen,Sarbach}
in the notation of our paper. The continuum identity
\begin{equation}
\psi'+{p\over r}\psi = {(r^p\psi)'\over r^p}
\end{equation}
suggests the finite differencing operator
\begin{equation}
(\tilde D\Psi)_i={(i+1)^p\Psi_{i+1}-(i-1)^p\Psi_{i-1}\over 2 i^p}.
\end{equation}
In \cite{CalabreseNeilsen,Sarbach} this is used on the interior points of a
centred grid. At the symmetry boundary
\begin{equation}
\label{Sarbachorigin}
(\tilde D\Psi)_0=(p+1)\Psi_i,
\end{equation}
and at the outer boundary 
\begin{equation}
(\tilde D\Psi)_M={M^p\Psi_M-(M-1)^p\Psi_{M-1}\over M^p}.
\end{equation}
This fits into our general approach with 
\begin{equation}
v_i=w_i=i^p \quad\hbox{for}\quad i\ne 0,
\end{equation}
$w_0=2/(1+p)$, and the outer boundary treatment
(\ref{WWtilde1},\ref{Dtildetmp1}), and hence is SBP.

It also appears to be second-order accurate, but it is not uniformly
so, in contrast to the methods derived here. As an example, for
$\psi=r$ (i.e. generic behaviour at the origin) and $p=2$, the local
error of the finite differencing operator $h^{-1}\tilde D$ is exactly
$h^2/r^2=1/i^2$. (For higher $p>0$, terms up to $(h/r)^p$ also
appear.) This does not go to zero with $h$ at fixed $i$. However, for
$p\ne 2,3$ the method converges with $h^2$ in the energy norm
\cite{Sarbachpc}.

%%%%%%%%%%%%%%%%%%%%%%%%%%%%%%%%%%%%%%%%%%%%%%%%%%%%%%%%%%%%%%%%%%%%%%%%%%

\section{Continuum boundary conditions involving derivatives}
\label{appendix:BCderivslwave}

Consider the class of boundary conditions of the form
\begin{equation}
\label{genBC}
\rho\pi+\sigma\psi+\mu\pi'+\nu\left(\psi'+{p\over r}\psi\right)=0,
\quad r=R
\end{equation} or equivalently
\begin{equation}
\label{genBCdot}
\rho\pi+\sigma\psi+\mu\dot\psi+\nu\dot\pi=0,
\quad r=R
\end{equation}
for $\rho,\sigma,\mu,\nu$ not all vanishing at once. To fix an overall
sign, we also assume that at least one of them is positive. We now use
an energy argument to show that these boundary conditions give rise to
a stable initial-boundary value problem if $\rho,\sigma,\mu,\nu\ge 0$
with $\rho\mu+\sigma\nu>0$. [The maximally dissipative special case
$\mu=\nu=0$ with $\rho\sigma\ge 0$ is also stable based on the energy
$E$ defined in (\ref{myE})].

We consider the energy
\begin{equation}
\label{Eb}
E_b\equiv E +{R^p\over 2s}(\mu\psi+\nu\pi)^2_{r=R},
\end{equation}
where $E$ is given by (\ref{myE}),
$E_b$ stands for $E$ modified by a boundary term, and $s$
is 
\begin{equation}
\label{outer_s}
s\equiv\rho\mu+\sigma\nu.  
\end{equation}
Its time
derivative is
\begin{eqnarray}
\label{Ebdot}
{dE_b\over dt}&=&R^p\left[\pi\psi
+{1\over s}(\mu\psi+\nu\pi)\left(\mu\dot\psi+\nu\dot\pi\right)
\right]_{r=R}
\nonumber \\
&=&R^p\left[\pi\psi
-{1\over s}(\mu\psi+\nu\pi)\left(\rho\pi+\sigma\psi\right)
\right]_{r=R}
\nonumber \\
&=&-{R^p\over s}(\mu\sigma\psi^2+\nu\rho\pi^2)_{r=R},
\end{eqnarray}
The necessary and
sufficient conditions for $E_b$ to be positive definite and its time
derivative to be non-positive are
\begin{equation}
\label{outersigns}
\rho\ge 0, \quad \sigma\ge 0, \quad
\mu\ge 0, \quad \nu \ge 0.
\end{equation}
We have $dE_b/dt=0$ if $\mu\sigma=\rho\pi=0$, and $dE_b/dt<0$
otherwise. However, the limiting case $\rho\mu=\sigma\nu=0$ is not
allowed because it would give $s=0$, except for the maximally
dissipative sub-case $\mu=\nu=0$, where $E$ and $dE/dt$ are given by
(\ref{myE}) and (\ref{myEdot}) instead of (\ref{Eb}) and
(\ref{Ebdot}).

%%%%%%%%%%%%%%%%%%%%%%%%%%%%%%%%%%%%%%%%%%%%%%%%%%%%%%%%%%%%%%%%%%%%%%

\section{Numerical boundary conditions involving derivatives}
\label{appendix:numericalBC} 

We define the modified numerical energy
\begin{equation}
\label{hatEb}
\hat E_b=\hat E + {\chi h^pM^p\over 2s}\left(\mu\Psi_M+\nu\Pi_M\right)^2,
\end{equation}
where $\chi$ parameterises finite differencing error in the
boundary term, as defined by (\ref{Bdef}). We find
\begin{eqnarray}
\label{dhatEbdt}
{d\hat E_b\over dt}&=&\chi h^pM^p\Bigl[\Pi_M\Psi_M 
\nonumber \\
&&+{1\over s}\left(\mu\Psi_M+\nu\Pi_M\right)
\left(\mu\dot\Psi_M+\nu\dot\Pi_M\right)\Bigr],
\end{eqnarray}
so if the numerical boundary could be chosen to be 
\begin{equation}
\rho\pi_M+\sigma\psi_M+\mu\dot\psi_M+\nu\dot\pi_M=0,
\end{equation}
the argument could be completed as in the continuum case. However, in
the notation of Appendix~\ref{appendix:Olsson}, $\dot
u={\cal PD}u$ and not ${\cal D}u$. We have not been able to find an
ansatz for $\cal L$ and $\hat E_b$ such that $d\hat E_b/dt\le0$.

Consider however the two subclasses of boundary conditions where
$dE_b/dt=0$ in the continuum. Consider first the case $\sigma=\nu=0$
with $\rho\mu>0$. Then (\ref{dhatEbdt}) reduces to
\begin{eqnarray}
{d\hat E_b\over dt}&=&\chi h^pM^p\left(\Pi_M\Psi_M 
+{\mu\over\rho}\Psi_M\dot\Psi_M\right) \nonumber \\
&=&\chi h^pM^p\left(\Pi_M\Psi_M 
+{\mu\over\rho}\Psi_Mh^{-1}(D\Pi)_M\right) \nonumber \\
&=& 0 ,
\end{eqnarray}
where the second equality holds because in this special case the
boundary condition is independent of $\Psi$ so that $\cal P$ only
acts on $\Pi$, and the last equality holds if we implement ${\cal
  L}u=0$ as
\begin{equation}
\rho\Pi_M+\mu h^{-1}(D\Pi)_M=0.
\end{equation}
The case $\rho=\mu=0$ with $\nu\sigma>0$ works the same way, with the
roles of $\Pi$ and $\Psi$ interchanged.

%%%%%%%%%%%%%%%%%%%%%%%%%%%%%%%%%%%%%%%%%%%%%%%%%%%%%%%%%%%%%%%%%%%%%%%%%%

\section{The projection method for imposing boundary conditions}
\label{appendix:Olsson}

For completeness, this Appendix summarises relevant methods from
\cite{Olsson}. Suppose a first-order in space and time system of PDEs
in one spatial dimension has been discretised in space as
\begin{equation}
\label{simpleupdate}
\dot u={\cal D}u
\end{equation}
Note that the vector $u$ in general ranges over multiple variables
(for example $\pi$ and $\psi$) as well as grid points (for example
$i$), and we use calligraphic letters such as $\cal D$ for operators
on this vector space.

Suppose this system has a discrete energy
\begin{equation}
{\hat E}\equiv{1\over 2}u^t{\cal W}u
\end{equation}
and obeys the SBP property that
\begin{equation}
{\cal B}\equiv {1\over 2}\left({\cal W}{\cal D}+{\cal D}^t{\cal W}\right)
\end{equation}
is a boundary operator. Then
\begin{equation}
{d{\hat E}\over dt}=u^t{\cal B}u
\end{equation}
is a boundary term.

We want to impose one or several
homogenous linear boundary conditions that we write as
\begin{equation}
{\cal L}u=0.
\end{equation}
In matrix notation where $u$ is a column vector, ${\cal L}$ is a
matrix that has one row for each boundary condition.

We define the inner product
\begin{equation}
(u,v)\equiv u^t{\cal W}v.
\end{equation}
In this notation we can write
\begin{equation}
\hat E={1\over 2}(u,u), \qquad {d\hat E\over dt}=(u,{\cal D}u).
\end{equation}
The adjoint with respect
to this inner product is defined by
\begin{equation}
({\cal A} u,v)\equiv(u,{\cal A}^\dagger v), 
\end{equation}
and is therefore given by
\begin{equation}
{\cal A}^\dagger={\cal W}^{-1}{\cal A}^t{\cal W}. 
\end{equation}

The operator
\begin{equation}
{\cal P}\equiv 1
-{\cal W}^{-1}{\cal L}^t({\cal L}{\cal W}^{-1}{\cal L}^t)^{-1}{\cal L}
\end{equation}
clearly obeys 
\begin{equation}
{\cal P}^2={\cal P}, \quad {\cal L}{\cal P}=0, \quad {\cal P}^\dagger={\cal P}, 
\end{equation}
and so is a self-adjoint projection operator into the space of grid
functions that obey the boundary conditions. If we now use the
semi-discrete evolution equation 
\begin{equation}
\dot u={\cal P}{\cal D}u
\end{equation}
instead of (\ref{simpleupdate}),
we have ${\cal L}\dot u=0$ exactly, and hence ${\cal L}u=0$ and
therefore ${\cal P}u=u$ at all
times if it holds initially. Then we have 
\begin{equation}
{d{\hat E}\over dt}=(u,{\cal P}{\cal D}u)
=({\cal P}u,{\cal D}u)
=(u,{\cal D}u)
\end{equation}
as before and so both the discrete energy bound and the desired
boundary conditions hold.

%%%%%%%%%%%%%%%%%%%%%%%%%%%%%%%%%%%%%%%%%%%%%%%%%%%%%%%%%%%%%%%%%%%%%%%%%%

%%%%%%%%%%%%%%%%%%%%%%%%%%%%%%%%%%%%%%%%%%%%%%%%%%%%%%%%%%%%%
%%%%%%%%%%%%%%%%%%%%%%%%%%%%%%%%%%%%%%%%%%%%%%%%%%%%%%%%%%%%%

%%%%%%%%%%%%%%%%%%%%%%%%%%%%%%%%%%%%%%%%%%%%%%%%%%%%%%%%%%%%
\begin{figure*}[ht!]
\psfrag{ipoint}{\large $i$}
\psfrag{d1}{\large $\!\!\!\delta_0$}
\psfrag{d2}{\large $\!\!\!\delta_1$}
\psfrag{d3}{\large $\!\!\!\delta_2$}
\hspace{-4mm}
\includegraphics[width=18cm]{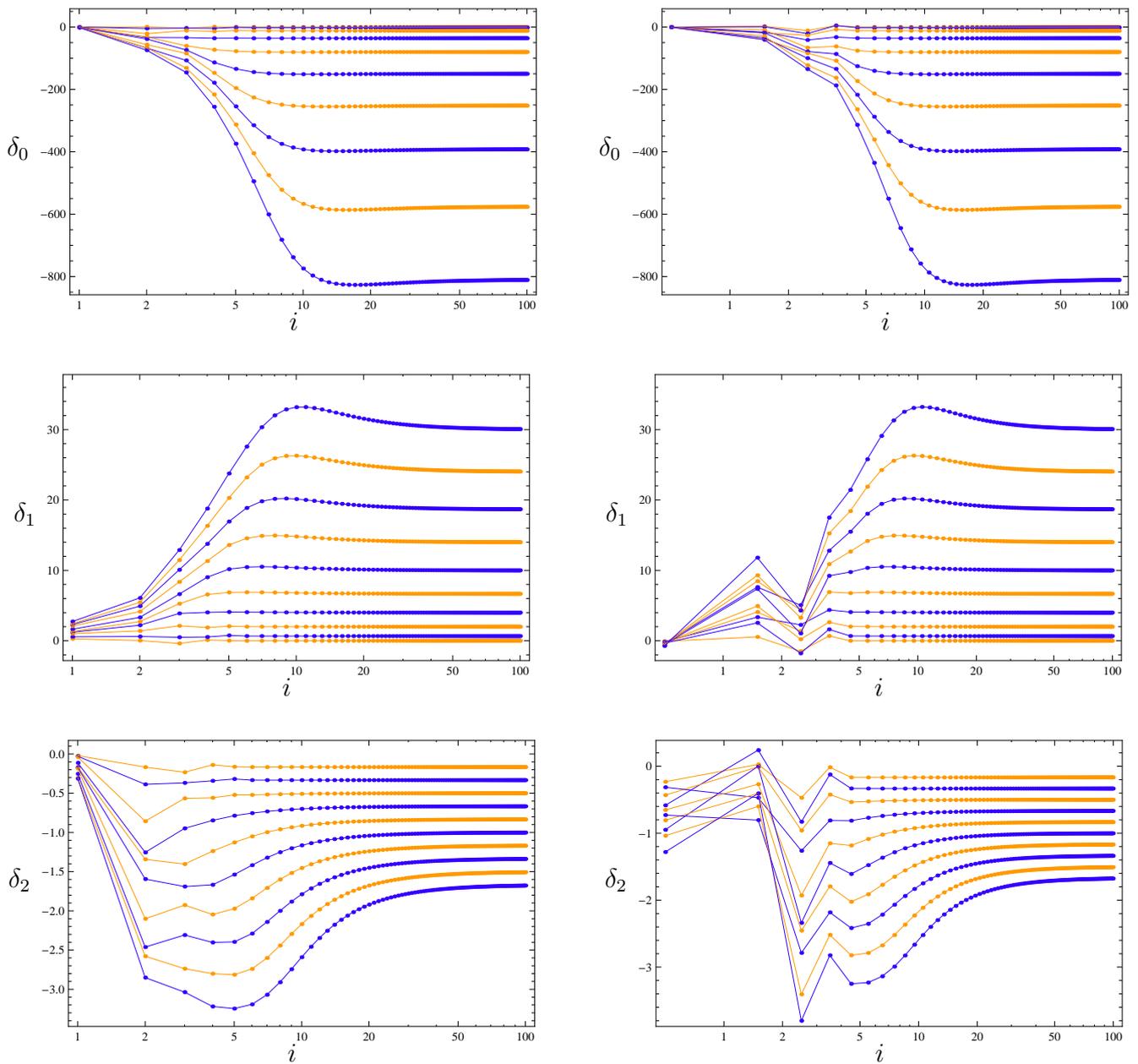}
\caption{\label{deltas_N=2} Values of $\delta_0$, $\delta_1$ and
  $\delta_2$ for SBP4 on the centred grid (left column) and on
  the staggered grid (right column), with $p=1,\ldots,10$. In all
  cases increasing values of $p$ correspond to lines further from the
  axis $\delta_i=0$, with even values of $p$ shown in blue (dark) and
  odd values in orange (light). Again, we see a rapid convergence
  towards their respective asymptotic values
  (\ref{asymp1b}--\ref{asymp3b}). The fact that $u_1$ (for the
  staggered grid grid) and $u_{3/2}$ or $u_{5/2}$ (for the centred
  grid) appear explicitly in the recurrence for a few low $i$ points,
  but not beyond, produces some irregular behaviour at those points.
}
\end{figure*}

%%%%%%%%%%%%%%%%%%%%%%%%%%%%%%%%%%%%%%%%%%%%%%%%%%%%%%%%%%%%

\end{document}